\documentclass[reqno,fleqn]{amsart}
\usepackage{amsfonts,amsthm,amssymb}
\usepackage{lipsum}
\usepackage[utf8]{inputenc}
\usepackage[english]{babel}

\usepackage{lmodern,textcomp}
\usepackage{amsmath}
\usepackage{mathtools}
\usepackage{latexsym}
\usepackage{tikz}
\usepackage{fancyvrb}
\usepackage{epsfig}
\usepackage{pstricks,slashbox,multirow}
\usepackage{graphicx}
\usepackage{rotating}
\usetikzlibrary{positioning}
\usepackage{subcaption}
\usepackage{datetime}
\newtheorem{theorem}{Theorem}

\newtheorem{definition}[theorem]{Definition}

\newtheorem{prm}{Problem}

\newtheorem{rem}[theorem]{Remark}

\title[A study on Type-2 isomorphic $C_n(R)$: Part 3: 384 pairs of Type-2 isomorphic $C_{32}(R)$]{A study on Type-2 isomorphic circulant graphs. \\ Part 3: 384 pairs of Type-2 isomorphic circulant graphs $C_{32}(R)$}

\author{\sc Vilfred Kamalappan} 
\address{Department of Mathematics, Central University of Kerala, Periye, Kasaragod, Kerala, India - 671 316.}
\email{vilfredkamal@gmail.com}



\date{}

\begin{document}

\begin{abstract} This article is the $3^{rd}$ part of a detailed study on Type-2 isomorphic circulant graphs having ten parts \cite{v2-1}-\cite{v2-10} and is a continuation of Part 2. Here, we obtain all 384 pairs of Type-2 isomorphic circulant graphs of order 32. 
\end{abstract}

\maketitle

	
\section{Introduction}

In \cite{v2-2}, we studied Type-2 isomorphic circulant graphs of orders 16, 24 and 27 and shown that the total number of pairs of Type-2 isomorphic circulant graphs of orders 16 and 24 are 8 and 32, respectively and the total number of triples of Type-2 isomorphic circulant graphs of order 27 are 12.  This paper is a continuation of \cite{v2-2} and using modified definition \ref{d3}, we study Type-2 isomorphic circulant graphs of order 32 and show that the total number of pairs of Type-2 isomorphic circulant graphs of orders 32 are 384. For basic definitions and results on isomorphic circulant graphs, refer \cite{v2-1, v2-2}.

\begin{definition}{\rm\cite{ad67}} \quad \label{d1} For $R =$ $\{r_1$, $r_2$, $\ldots$, $r_k\}$ and $S$ = $\{s_1$, $s_2$, $\ldots$, $s_k\}$, circulant graphs $C_n(R)$ and $C_n(S)$ are {\it Adam's isomorphic} if there exists a positive integer $x$ $\ni$ $\gcd(n, x)$ = 1 and $S$ = $\{xr_1$, $xr_2$, $\ldots$, $xr_k\}_n^*$ where $<r_i>_n^*$, the {\it reflexive modular reduction} of a sequence $< r_i >$, is the sequence obtained by reducing each $r_i$ under modulo $n$ to yield $r_i'$ and then replacing all resulting terms $r_i'$ which are larger than $\frac{n}{2}$ by $n-r_i'.$  
\end{definition}

A circulant graph $C_n(R)$ is said to have {\em Cayley Isomorphism (CI) property} if whenever $C_n(S)$ is isomorphic to $C_n(R)$, they are Adam’s isomorphic \cite{krsi}.

\begin{theorem} \cite{v24} \label{t2} Let $Ad_n(C_n(R))$ = $\{\varphi_{n,x}(C_n(R)) = C_n(xR): x\in\varphi_n \}$. Then, $C_n(S)\in Ad_n(C_n(R))$ if and only if $Ad_n(C_n(R))$ = $Ad_n(C_n(S))$ if and only if $C_n(R)\in Ad_n(C_n(S))$. \hfill $\Box$
\end{theorem}

In \cite{v2-1}, Vilfred modified the definition of Type-2 isomorphism of $C_n(R)$ w.r.t. $m$  given in \cite{v2-2-arX} as follows and hereafter we use the definition \ref{d3}. 

\begin{definition} \cite{v2-1} \quad  \label{d3} Let $V(K_n) = \{u_0,u_1,u_2,...,u_{n-1}\}$, $V(C_n(R))$ = $\{v_0, v_1, v_2, ... , v_{n-1}\}$, $|R| \geq 3$, $r\in R$ and $m > 1$ and $m^3$ be divisors of $\gcd(n, r)$ and $n$, respectively. Define 1-1 mapping $\theta_{n,m,t} :$ $V(C_n(R)) \rightarrow V(K_n)$ such that $\theta_{n,m,t}(v_x) = u_{x+jtm}$,  $\theta_{n,m,t}((v_x, v_{x+s}))$ = $(\theta_{n,m,t}(v_x),$ $\theta_{n,m,t}(v_{x+s}))$ under subscript arithmetic modulo $n$ and $\theta_{n,m,t}(C_n(R))$ = $C_n(\theta_{n,m,t}(R))$ where $\theta_{n,m,t}(R)$ in $C_n(\theta_{n,m,t}(R))$ is calculated under the reflexive modulo $n$, $\forall$ $x \in \mathbb{Z}_n$, $x = qm+j,$ $0 \leq j \leq m-1$, $s\in R$ and $0 \leq q,t \leq \frac{n}{m} -1$. And for a particular value of $t,$ if  $\theta_{n,m,t}(C_n(R))$ = $C_n(S)$ for some $S$  and  $S \neq yR$ for all $y\in \varphi_n$ under reflexive modulo $n,$ then $C_n(R)$ and $C_n(S)$ are called {\em isomorphic circulant graphs of Type-2 w.r.t. $m$.} 
	
When $C_n(R)$ and $C_n(S)$ are Type-2 isomorphic w.r.t. $m$, then we also say that $C_{kn}(kR)$  and $C_{kn}(kS)$ are Type-2 isomorphic w.r.t. $m$, $k\in\mathbb{N}$. Here, $k.C_n(T)$ = $C_{kn}(kT)$, $k\in\mathbb{N}$. 	 
\end{definition}

\begin{rem} \cite{v2-1}  \label{r4} Following steps are used to establish Type-2 isomorphism w.r.t. $m$ between circulant graphs $C_n(R)$ and $C_n(S)$. (i) $R$ $\neq$ $S$ and $|R| = |S| \geq 3$; (ii) $\exists$ $r\in R,S$ and $m > 1$ $\ni$ $m$ is a divisor of $\gcd(n, r)$, $m^3$ is a divisor of $n$ and for some $t$ $\ni$ $1 \leq t \leq \frac{n}{m} -1$, $\theta_{n,m,t}(C_n(R))$ = $C_n(S)$ and (iii) $S$ $\neq$ $xR$ for all $x\in\varphi_n$ under arithmetic reflexive modulo $n$. 

Thus, if $C_n(R)$ and $C_n(S)$ are Type-2 isomorphic circulant graphs w.r.t. $m$, then there exist $r\in R,S$, $m > 1$ and some $t$ $\ni$ $m$ is a divisor of $\gcd(n, r)$, $m^3$ is a divisor of $n$, $1 \leq t \leq \frac{n}{m} -1$, $\theta_{n,m,t}(C_n(R))$ = $C_n(S)$ and $S$ $\neq$ $xR$ for all $x\in\varphi_n$ under arithmetic reflexive modulo $n$.
\end{rem} 

\begin{rem}  \label{r5} \quad The calculation on jump sizes $r_i$s which are integer multiples of $m$ need not be done under the transformation $\theta_{n,m,t}$, while searching for possible value(s) of $t$ for which the transformed graph $\theta_{n,m,t}(C_n(R))$ is circulant of the form $C_n(S)$ for some $S \subseteq [1, \frac{n}{2}]$, as there is no change in these $r_i$s where $r\in R$ and $m > 1$ and $m^3$ are divisors of $\gcd(n, r)$ and $n$, respectively. 

Thus, if $\theta_{n,m,t}(C_n(R))$ = $C_n(S)$ for some $S$ and thereby $C_n(R)$ $\cong$ $C_n(S)$, then $\theta_{n,m,t}(C_n(R \cup mT))$ = $C_n(S \cup mT)$ for any $T$ and thereby $C_n(R \cup mT)$ $\cong$ $C_n(S \cup mT)$.

Also, for a given $C_n(R)$, w.r.t. different values of $m$, we may get different Type-2 isomorphic circulant graphs.
\end{rem}

\begin{rem} \cite{v2-1} \label{r6} \quad {\rm For given $C_n(R)$ and $C_n(S)$ when either $\theta_{n,m,t}(C_n(R))$ = $C_n(S)$ for some $t$ or $C_n(xR)$ = $C_n(S)$ for some $x$, then $C_n(R)$ and $C_n(S)$ are isomorphic, $0 \leq t \leq \frac{n}{m} -1$ and $x\in\varphi_n$. }
\end{rem}

\begin{theorem} \cite{v2-6} \label{t7} {\rm Let $p$ be an odd prime number, $1 \leq i \leq p$, $1 \leq x \leq p-1$, $y\in\mathbb{N}_0$, $0 \leq y \leq np-1$, $1 \leq x+yp \leq np^2-1$, $d^{np^3, x+yp}_i = (i-1)xpn+x+yp$,  $R^{np^3, x+yp}_i$ $=$ $\{p$, $d^{np^3, x+yp}_i$, $np^2-d^{np^3, x+yp}_i$, $np^2+d^{np^3, x+yp}_i$, $2np^2-d^{np^3, x+yp}_i$, $2np^2+$ $d^{np^3, x+yp}_i,$ $3np^2-d^{np^3, x+yp}_i$, $3np^2+d^{np^3, x+yp}_i$, . . . , $(p-1)np^2$ - $d^{np^3, x+yp}_i$, $(p-1)np^2+d^{np^3, x+yp}_i$, $np^3-d^{np^3, x+yp}_i$, $np^3-p\}$ and $i,j,n,x\in\mathbb{N}$. Then, for a given set of values of $n$, $p$, $x$ and $y$, $\theta_{np^3,p,jn} (C_{np^3}(R^{np^3, x+yp}_i))$ = $C_{np^3}(R^{np^3, x+yp}_{i+j})$ and the $p$ circulant graphs $C_{np^3}(R^{np^3, x+yp}_i)$ are isomorphic of Type-2 w.r.t.  $p$, $1 \leq i,j \leq p$ where $i+j$ in $R^{np^3, x+yp}_{i+j}$ is calculated under addition modulo $p$ and $C_{np^3}(R^{np^3, x+yp}_0)$ = $C_{np^3}(R^{np^3, x+yp}_p)$. \hfill $\Box$}
\end{theorem}

\section{Main results}

Given a circulant graph $C_n(R)$ having isomorphic circulant graphs of Type-2 w.r.t. $m$, remark \ref{r5} helps us to obtain more  isomorphic graphs which covers Type-2 w.r.t. $m$ as well as some Type-1 isomorphic graphs of $C_n(R)$. In this section, we study Type-2 isomorphic circulant graphs of order 32. 

\subsection{On Type-2 isomorphic circulant graphs of $C_{32}(r_1,r_2,r_3)$ w.r.t. $m$ = 2}

In this subsection, we obtain all pairs of Type-2 isomorphic circulant graphs of order 32. We show that the total number of pairs of Type-2 isomorphic circulant graphs of order 32 are 384. In the next two problems, we present our findings on Type-2 isomorphic circulant graphs of $C_{32}(R)$ and show that the following 16 pairs of circulant graphs of the form $C_{32}(r_1,r_2,r_3)$ are isomorphic and only 8 pairs among them are of Type-2 isomorphic w.r.t. $m$ = 2 where 32 = $4\times 2^3$.  
 
\begin{prm}\quad \label{p3.5} {\rm Show that the following statements are true. 
\begin{enumerate}
\item [\rm (1a)]   $C_{32}(1,2,15)$ and $C_{32}(2,7,9)$ are  Type-2 isomorphic w.r.t. $m$ = 2; 
\item [\rm (1b)]   $C_{32}(1,4,15)$ and $C_{32}(4,7,9)$ are  Type-1 isomorphic;
\item [\rm (1c)]   $C_{32}(1,6,15)$ and $C_{32}(6,7,9)$ are  Type-2 isomorphic w.r.t. $m$ = 2;
\item [\rm (1d)]   $C_{32}(1,8,15)$ and $C_{32}(7,8,9)$ are  Type-1 isomorphic;
\item [\rm (1e)]   $C_{32}(1,10,15)$ and $C_{32}(7,9,10)$ are  Type-2 isomorphic w.r.t. $m$ = 2;
\item [\rm (1f)]   $C_{32}(1,12,15)$ and $C_{32}(7,9,12)$ are  Type-1 isomorphic;
\item [\rm (1g)]   $C_{32}(1,14,15)$ and $C_{32}(7,9,14)$ are  Type-2 isomorphic w.r.t. $m$ = 2;
\item [\rm (1h)]   $C_{32}(1,15,16)$ and $C_{32}(7,9,16)$ are  Type-1 isomorphic;
\item [\rm (2a)]   $C_{32}(2,3,13)$ and $C_{32}(2,5,11)$ are  Type-2 isomorphic w.r.t. $m$ = 2;
\item [\rm (2b)]   $C_{32}(3,4,13)$ and $C_{32}(4,5,11)$ are  Type-1 isomorphic;
\item [\rm (2c)]   $C_{32}(3,6,13)$ and $C_{32}(5,6,11)$ are  Type-2 isomorphic w.r.t. $m$ = 2;
\item [\rm (2d)]   $C_{32}(3,8,13)$ and $C_{32}(5,8,11)$ are  Type-1 isomorphic;
\item [\rm (2e)]   $C_{32}(3,10,13)$ and $C_{32}(5,10,11)$ are  Type-2 isomorphic w.r.t. $m$ = 2;
\item [\rm (2f)]   $C_{32}(3,12,13)$ and $C_{32}(5,11,12)$ are  Type-1 isomorphic;
\item [\rm (2g)]   $C_{32}(3,13,14)$ and $C_{32}(5,11,14)$ are  Type-2 isomorphic w.r.t. $m$ = 2;
\item [\rm (2h)]   $C_{32}(3,13,16)$ and $C_{32}(5,11,16)$ are  Type-1 isomorphic.
\end{enumerate}
  }
\end{prm}
\noindent
{\bf Solution.}  Here, $n$ = 32 = $4\times 2^3$ and so the possible values of $m > 1$ $\ni$ $m$ is a divisor of $\gcd(n, r)$ = $\gcd(32, r)$, $m^3$ is a divisor of 32 and $r\in R$ for the existence of isomorphic circulant graphs $C_{32}(R)$ of Type-2 w.r.t. $m$ is $m$ = 2 and 2 = $\gcd(32, 2)$ = $\gcd(32, 6)$ = $\gcd(32, 10)$ = $\gcd(32, 14)$, 4 = $\gcd(32, 4)$ = $\gcd(32, 12)$ and 8 = $\gcd(32, 8)$. Also, we have $\theta_{32,2,t}(s)$ = $s+2jt$ where $s\in\mathbb{Z}_{32}$, $s$ = $qm+j$, $j,m,q\in\mathbb{N}_0$, $m$ = 2, $0 \leq j \leq 1$ and $0 \leq t \leq \frac{32}{\gcd(32, m)}-1$ = $\frac{32}{\gcd(32, 2)}-1$ = 15.
\begin{enumerate}
\item [\rm (1a)]   $\theta_{32,2, 4}(C_{32} (1,2,15))$ = $\theta_{32,2, 4}(C_{32} (1,2,15,  17,30,31))$  = $C_{32} (\theta_{32,2, 4}(1,2,15,  17,30,31))$ 

\hspace{2.5cm}  = $C_{32} (9,2,23, 25,30,7)$ = $C_{32} (2,7,9,  23,25,30)$  

\hspace{2.5cm} = $C_{32} (2,7,9)$.  $\Rightarrow$ $C_{32} (1,2,15)$ $\cong$ $C_{32} (2,7,9)$.
\\
$Ad_{32} (C_{32} (1,2,15))$ = $\{\varphi_{32,x}(C_{32} (1,2,15): x = 1,3,5,7,9,11,13,15\}$ 

\hspace{2cm} = $\{C_{32} (x(1,2,15)): x = 1,3,5,7,9,11,13,15\}$ 

\hspace{2cm}  = $\{C_{32} (1,2,15), C_{32} (3,6,13), C_{32}(5,10,11), C_{32}(7,9,14)\}$ 

\hspace{2cm} = $\{C_{32} (x(1,2,15)): x = 1,3,5,7\}$. 

This implies, $C_{32} (2,7,9) \notin Ad_{32} (C_{32} (1,2,15))$. 

Hence, $C_{32} (1,2,15)$ and $C_{32} (2,7,9)$ are Type-2 isomorphic w.r.t. $m$ = 2.  

\item [\rm (1b)]   $C_{32} (7(1,4,15))$ = $C_{32} (4,7,9)$.  $\Rightarrow$ $C_{32} (1,4,15)$ and $C_{32} (4,7,9)$ are Type-1 isomorphic.  

\item [\rm (1c)]   $\theta_{32,2, 4}(C_{32} (1,6,15))$ = $C_{32} (6,7,9)$.  $\Rightarrow$ $C_{32} (1,6,15)$ $\cong$ $C_{32} (6,7,9)$.
\\
$Ad_{32} (C_{32} (1,6,15))$ = $\{C_{32} (1,6,15), C_{32} (3,13,14), C_{32}(2,5,11), C_{32}(7,9,10)\}$.

This implies, $C_{32} (6,7,9) \notin Ad_{32} (C_{32} (1,6,15))$. 

Hence, $C_{32} (1,6,15)$ and $C_{32} (6,7,9)$ are Type-2 isomorphic w.r.t. $m$ = 2.  

\item [\rm (1d)]   $C_{32} (7(1,8,15))$ = $C_{32} (7,8,9)$.  $\Rightarrow$ $C_{32} (1,8,15)$ and $C_{32} (7,8,9)$ are Type-1 isomorphic.  

\item [\rm (1e)]   $\theta_{32,2, 4}(C_{32} (1,10,15))$ = $C_{32} (7,9,10)$.  $\Rightarrow$ $C_{32} (1,10,15)$ $\cong$ $C_{32} (7,9,10)$.
\\
$Ad_{32} (C_{32} (1,10,15))$ = $\{C_{32} (1,10,15), C_{32} (2,3,13), C_{32}(5,11,14), C_{32}(6,7,9)\}$.

This implies, $C_{32} (7,9,10) \notin Ad_{32} (C_{32} (1,10,15))$. 

Hence, $C_{32} (1,10,15)$ and $C_{32} (7,9,10)$ are Type-2 isomorphic w.r.t. $m$ = 2.  

\item [\rm (1f)]   $C_{32} (7(1,12,15))$ = $C_{32} (7,9,12)$.  $\Rightarrow$ $C_{32} (1,12,15)$ and $C_{32} (7,9,12)$ are Type-1 isomorphic.  

\item [\rm (1g)]   $\theta_{32,2, 4}(C_{32} (1,14,15))$ = $C_{32} (7,9,14)$.  $\Rightarrow$ $C_{32} (1,14,15)$ $\cong$ $C_{32} (7,9,14)$.
\\
$Ad_{32} (C_{32} (1,14,15))$ = $\{C_{32} (1,14,15), C_{32} (3,10,13), C_{32}(5,6,11), C_{32}(2,7,9)\}$.

This implies, $C_{32} (7,9,14) \notin Ad_{32} (C_{32} (1,14,15))$. 

Hence, $C_{32} (1,14,15)$ and $C_{32} (7,9,14)$ are Type-2 isomorphic w.r.t. $m$ = 2.  

\item [\rm (1h)]   $C_{32} (7(1,15,16))$ = $C_{32} (7,9,16)$.  $\Rightarrow$ $C_{32} (1,15,16)$ and $C_{32} (7,9,16)$ are Type-1 isomorphic.  

\item [\rm (2a)]   $\theta_{32,2, 4}(C_{32} (2,3,13))$ = $C_{32} (2,5,11)$.  $\Rightarrow$ $C_{32} (2,3,13)$ $\cong$ $C_{32} (2,5,11)$.
\\
$Ad_{32} (C_{32} (2,3,13))$ = $\{C_{32} (2,3,13), C_{32} (6,7,9), C_{32}(1,10,15), C_{32}(5,11,14)\}$.

This implies, $C_{32} (2,5,11) \notin Ad_{32} (C_{32} (2,3,13))$. 

Hence, $C_{32} (2,3,13)$ and $C_{32} (2,5,11)$ are Type-2 isomorphic w.r.t. $m$ = 2.  

\item [\rm (2b)]   $C_{32} (7(3,4,13))$ = $C_{32} (4,5,11)$.  $\Rightarrow$ $C_{32} (3,4,13)$ and $C_{32} (4,5,11)$ are Type-1 isomorphic.  

\item [\rm (2c)]   $\theta_{32,2, 4}(C_{32} (3,6,13))$ = $C_{32} (5,6,11)$.  $\Rightarrow$ $C_{32} (3,6,13)$ $\cong$ $C_{32} (5,6,11)$.
\\
$Ad_{32} (C_{32} (3,6,13))$ = $\{C_{32} (3,6,13), C_{32} (7,9,14), C_{32}(1,2,15), C_{32}(5,10,11)\}$.

This implies, $C_{32} (5,6,11) \notin Ad_{32} (C_{32} (3,6,13))$. 

Hence, $C_{32} (3,6,13)$ and $C_{32} (5,6,11)$ are Type-2 isomorphic w.r.t. $m$ = 2.  

\item [\rm (2d)]   $C_{32} (7(3,8,13))$ = $C_{32} (5,8,11)$.  $\Rightarrow$ $C_{32} (3,8,13)$ and $C_{32} (5,8,11)$ are Type-1 isomorphic.  

\item [\rm (2e)]   $\theta_{32,2, 4}(C_{32} (3,10,13))$ = $C_{32} (5,10,11)$.  $\Rightarrow$ $C_{32} (3,10,13)$ $\cong$ $C_{32} (5,10,11)$.
\\
$Ad_{32} (C_{32} (3,10,13))$ = $\{C_{32} (3,10,13), C_{32} (2,7,9), C_{32}(1,14,15), C_{32}(5,6,11)\}$.

This implies, $C_{32} (5,10,11) \notin Ad_{32} (C_{32} (3,10,13))$. 

Hence, $C_{32} (3,10,13)$ and $C_{32} (5,10,11)$ are Type-2 isomorphic w.r.t. $m$ = 2.  

\item [\rm (2f)]   $C_{32} (7(3,12,13))$ = $C_{32} (5,11,12)$.  $\Rightarrow$ $C_{32} (3,12,13)$ and $C_{32} (5,11,12)$ are Type-1 isomorphic.  

\item [\rm (2g)]   $\theta_{32,2, 4}(C_{32} (3,13,14))$ = $C_{32} (5,11,14)$.  $\Rightarrow$ $C_{32} (3,13,14)$ $\cong$ $C_{32} (5,11,14)$.
\\
$Ad_{32} (C_{32} (3,13,14))$ = $\{C_{32} (3,13,14), C_{32} (7,9,10), C_{32}(1,6,15), C_{32}(2,5,11)\}$.

This implies, $C_{32} (5,11,14) \notin Ad_{32} (C_{32} (3,13,14))$. 

Hence, $C_{32} (3,13,14)$ and $C_{32} (5,11,14)$ are Type-2 isomorphic w.r.t. $m$ = 2.  

\item [\rm (2h)]   $C_{32} (7(3,13,16))$ = $C_{32} (5,11,16)$.  $\Rightarrow$ $C_{32} (3,13,16)$ and $C_{32} (5,11,16)$ are Type-1 isomorphic.  \hfill $\Box$ 
\end{enumerate}

\begin{prm}\quad \label{p2} {\rm Show that the following statements are true. 
\begin{enumerate}
\item [\rm (b1)]   $\theta_{32,2, 4}(C_{32}(1,4,15))$ = $C_{32}(4,7,9)$ and $C_{32}(1,4,15)$ and $C_{32}(4,7,9)$ are Type-1 isomorphic;
\item [\rm (b2)]   $\theta_{32,2, 4}(C_{32}(1,4,8,15))$ = $C_{32}(4,7,8,9)$ and

$C_{32}(1,4,8,15)$ and $C_{32}(4,7,8,9)$ are Type-1 isomorphic;
\item [\rm (b3)]   $\theta_{32,2, 4}(C_{32}(1,4,12,15))$ = $C_{32}(4,7,9,12)$ and 

$C_{32}(1,4,12,15)$ and $C_{32}(4,7,9,12)$ are Type-1 isomorphic;
\item [\rm (b4)]   $\theta_{32,2, 4}(C_{32}(1,4,15,16))$ = $C_{32}(4,7,9,16)$ and 

$C_{32}(1,4,15,16)$ and $C_{32}(4,7,9,16)$ are Type-1 isomorphic;
\item [\rm (b5)]   $\theta_{32,2, 4}(C_{32}(1,8,12,15))$ = $C_{32}(7,8,9,12)$ and 

$C_{32}(1,8,12,15)$ and $C_{32}(7,8,9,12)$ are Type-1 isomorphic;
\item [\rm (b6)]   $\theta_{32,2, 4}(C_{32}(1,12,15,16))$ = $C_{32}(7,9,12,16)$ and 

$C_{32}(1,12,15,16)$ and $C_{32}(7,9,12,16)$ are Type-1 isomorphic;
\item [\rm (b7)]   $\theta_{32,2, 4}(C_{32}(1,4,8,12,15))$ = $C_{32}(4,7,8,9,12)$ and 

$C_{32}(1,4,8,12,15)$ and $C_{32}(4,7,8,9,12)$ are Type-1 isomorphic;
\item [\rm (b8)]   $\theta_{32,2, 4}(C_{32}(1,4,8,15,16))$ = $C_{32}(4,7,8,9,16)$ and 

$C_{32}(1,4,8,15,16)$ and $C_{32}(4,7,8,9,16)$ are Type-1 isomorphic;
\item [\rm (b9)]   $\theta_{32,2, 4}(C_{32}(1,4,12,15,16))$ = $C_{32}(4,7,9,12,16)$ and 

$C_{32}(1,4,12,15,16)$ and $C_{32}(4,7,9,12,16)$ are Type-1 isomorphic;
\item [\rm (b10)]   $\theta_{32,2, 4}(C_{32}(1,8,12,15,16))$ = $C_{32}(7,8,9,12,16)$ and 

$C_{32}(1,8,12,15,16)$ and $C_{32}(7,8,9,12,16)$ are Type-1 isomorphic;
\item [\rm (b11)]   $\theta_{32,2, 4}(C_{32}(1,4,8,12,15,16))$ = $C_{32}(4,7,8,9,12,16)$ and 

$C_{32}(1,4,8,12,15,16)$ and $C_{32}(4,7,8,9,12,16)$ are Type-1 isomorphic;
\item [\rm (d1)]   $\theta_{32,2, 4}(C_{32}(1,8,15))$ = $C_{32}(7,8,9)$ and 

$C_{32}(1,8,15)$ and $C_{32}(7,8,9)$ are Type-1 isomorphic;
\item [\rm (d2)]   $\theta_{32,2, 4}(C_{32}(1,8,15,16))$ = $C_{32}(7,8,9,16)$ and 

$C_{32}(1,8,15,16)$ and $C_{32}(7,8,9,16)$ are Type-1 isomorphic;
\item [\rm (f1)]   $\theta_{32,2, 4}(C_{32}(1,12,15))$ = $C_{32}(7,9,12)$ and 

$C_{32}(1,12,15)$ and $C_{32}(7,9,12)$ are Type-1 isomorphic;
\item [\rm (h1)]   $\theta_{32,2, 4}(C_{32}(1,15,16))$ = $C_{32}(7,9,16)$ and 

$C_{32}(1,15,16)$ and $C_{32}(7,9,16)$ are Type-1 isomorphic;
\item [\rm (i1)]   $\theta_{32,2, 4}(C_{32}(3,4,13))$ = $C_{32}(4,5,11)$ and 

$C_{32}(3,4,13)$ and $C_{32}(4,5,11)$ are Type-1 isomorphic;
\item [\rm (i2)]   $\theta_{32,2, 4}(C_{32}(3,4,8,13))$ = $C_{32}(4,5,8,11)$ and 

$C_{32}(3,4,8,13)$ and $C_{32}(4,5,8,11)$ are Type-1 isomorphic;
\item [\rm (i3)]   $\theta_{32,2, 4}(C_{32}(3,4,12,13))$ = $C_{32}(4,5,11,12)$ and 

$C_{32}(3,4,12,13)$ and $C_{32}(4,5,11,12)$ are Type-1 isomorphic;
\item [\rm (i4)]   $\theta_{32,2, 4}(C_{32}(3,4,13,16))$ = $C_{32}(4,5,11,16)$ and 

$C_{32}(3,4,13,16)$ and $C_{32}(4,5,11,16)$ are Type-1 isomorphic;
\item [\rm (i5)]   $\theta_{32,2, 4}(C_{32} (3,8,12,13))$ = $C_{32} (5,8,11,12)$ and  

$C_{32} (3,8,12,13)$ and $C_{32} (5,8,11,12)$ are Type-1 isomorphic;
\item [\rm (i6)]   $\theta_{32,2, 4}(C_{32} (3,12,13,16))$ = $C_{32} (5,11,12,16)$ and  

$C_{32} (3,12,13,16)$ and $C_{32} (5,11,12,16)$ are Type-1 isomorphic;

\item [\rm (i7)]   $\theta_{32,2, 4}(C_{32}(3,4,8,12,13))$ = $C_{32}(4,5,8,11,12)$ and 

$C_{32}(3,4,8,12,13)$ and $C_{32}(4,5,8,11,12)$ are Type-1 isomorphic;
\item [\rm (i8)]   $\theta_{32,2, 4}(C_{32}(3,4,8,13,16))$ = $C_{32}(4,5,8,11,16)$ and 

$C_{32}(3,4,8,13,16)$ and $C_{32}(4,5,8,11,16)$ are Type-1 isomorphic;
\item [\rm (i9)]   $\theta_{32,2, 4}(C_{32}(3,4,12,13,16))$ = $C_{32}(4,5,11,12,16)$ and 

$C_{32}(3,4,12,13,16)$ and $C_{32}(4,5,11,12,16)$ are Type-1 isomorphic;
\item [\rm (i10)]   $\theta_{32,2, 4}(C_{32}(3,8,12,13,16))$ = $C_{32}(5,8,11,12,16)$ and 

$C_{32}(3,8,12,13,16)$ and $C_{32}(5,8,11,12,16)$ are Type-1 isomorphic;
\item [\rm (i11)]   $\theta_{32,2, 4}(C_{32}(3,4,8,12,13,16))$ = $C_{32}(4,5,8,11,12,16)$ and 

$C_{32}(3,4,8,12,13,16)$ and $C_{32}(4,5,8,11,12,16)$ are Type-1 isomorphic;
\item [\rm (j1)]   $\theta_{32,2, 4}(C_{32}(3,8,13))$ = $C_{32}(5,8,11)$ and $C_{32}(3,8,13)$ and $C_{32}(5,8,11)$ are Type-1 isomorphic;
\item [\rm (j2)]   $\theta_{32,2, 4}(C_{32}(3,8,13,16))$ = $C_{32}(5,8,11,16)$ and 

$C_{32}(3,8,13,16)$ and $C_{32}(5,8,11,16)$ are Type-1 isomorphic;
\item [\rm (k1)]  $\theta_{32,2, 4}(C_{32}(3,12,13))$ = $C_{32}(5,11,12)$ and 

$C_{32}(3,12,13)$ and $C_{32}(5,11,12)$ are Type-1 isomorphic;
\item [\rm (l1)]   $\theta_{32,2, 4}(C_{32}(3,13,16))$ = $C_{32}(5,11,16)$ and 

$C_{32}(3,13,16)$ and $C_{32}(5,11,16)$ and are Type-1 isomorphic.
\end{enumerate}
  }
\end{prm}
\noindent
{\bf Solution.} \quad In each case, it is easy to see that  $\theta_{32,2, 4}(C_{32}(R))$ = $C_{32}(S)$ and we will show that $C_{32}(S)$ = $C_{32}(7R)$,  $7\in\varphi_{32}$ and thereby $C_{32}(R)$ and $C_{32}(S)$ are Type-1 isomorphic.
\begin{enumerate}
\item [\rm (b1)]   $C_{32} (7(1,4,15))$ = $C_{32} (4,7,9)$.  $\Rightarrow$ $C_{32} (1,4,15)$ and $C_{32} (4,7,9)$ are Type-1 isomorphic.
\item [\rm (b2)]   $C_{32} (7(1,4,8,15))$ = $C_{32} (4,7,8,9)$.  $\Rightarrow$ $C_{32} (1,4,8,15)$ and $C_{32} (4,7,8,9)$ are Type-1 isomorphic.
\item [\rm (b3)]   $C_{32} (7(1,4,12,15))$ = $C_{32} (4,7,9,12)$.  
\\
$\Rightarrow$ $C_{32} (1,4,12,15)$ and $C_{32} (4,7,9,12)$ are Type-1 isomorphic.
\item [\rm (b4)]    $C_{32} (7(1,4,15,16))$ = $C_{32} (4,7,9,16)$.  
\\
$\Rightarrow$ $C_{32} (1,4,15,16)$ and $C_{32} (4,7,9,16)$ are Type-1 isomorphic. 
\item [\rm (b5)]   $C_{32}(7(1,8,12,15))$ = $C_{32}(7,8,9,12)$.  
\\
$\Rightarrow$ $C_{32}(1,8,12,15)$ and $C_{32}(7,8,9,12)$ are Type-1 isomorphic;
\item [\rm (b6)]   $C_{32}(7(1,12,15,16))$ = $C_{32}(7,9,12,16)$ and 
\\
$\Rightarrow$ $C_{32}(1,12,15,16)$ and $C_{32}(7,9,12,16)$ are Type-1 isomorphic;
\item [\rm (b7)]   $C_{32} (7(1,4,8,12,15))$ = $C_{32} (4,7,8,9,12)$.  
\\
$\Rightarrow$ $C_{32} (1,4,8,12,15)$ and $C_{32} (4,7,8,9,12)$ are Type-1 isomorphic.
\item [\rm (b8)]   $C_{32} (7(1,4,8,15,16))$ = $C_{32} (4,7,8,9,16)$.  
\\
$\Rightarrow$ $C_{32} (1,4,8,15,16)$ and $C_{32} (4,7,8,9,16)$ are Type-1 isomorphic.
\item [\rm (b9)]   $C_{32} (7(1,4,12,15,16))$ = $C_{32} (4,7,9,12,16)$.  
\\
$\Rightarrow$ $C_{32} (1,4,12,15,16)$ and $C_{32} (4,7,9,12,16)$ are Type-1 isomorphic.
\item [\rm (b10)]  $C_{32} (7(1,8,12,15,16))$ = $C_{32} (7,8,9,12,16)$.  
\\
$\Rightarrow$ $C_{32} (1,8,12,15,16)$ and $C_{32} (7,8,9,12,16)$ are Type-1 isomorphic. 
\item [\rm (b11)]  $C_{32} (7(1,4,8,12,15,16))$ = $C_{32} (4,7,8,9,12,16)$.  
\\
$\Rightarrow$ $C_{32} (1,4,8,12,15,16)$ and $C_{32} (4,7,8,9,12,16)$ are Type-1 isomorphic.
\item [\rm (d1)]  $C_{32} (7(1,8,15))$ = $C_{32} (7,8,9)$.  $\Rightarrow$ $C_{32} (1,8,15)$ and $C_{32} (7,8,9)$ are Type-1 isomorphic. 
\item [\rm (d2)]  $C_{32} (7(1,8,15,16))$ = $C_{32} (7,8,9,16)$.  
\\
$\Rightarrow$ $C_{32} (1,8,15,16)$ and $C_{32} (7,8,9,16)$ are Type-1 isomorphic. 
\item [\rm (f1)]  $C_{32} (7(1,12,15))$ = $C_{32} (7,9,12)$.  $\Rightarrow$ $C_{32} (1,12,15)$ and $C_{32} (7,9,12)$ are Type-1 isomorphic. 
\item [\rm (h1)]  $C_{32} (7(1,15,16))$ = $C_{32} (7,9,16)$.  $\Rightarrow$ $C_{32} (1,15,16)$ and $C_{32} (7,9,16)$ are Type-1 isomorphic. 
\item [\rm (i1)]  $C_{32} (7(3,4,13))$ = $C_{32} (4,5,11)$.  $\Rightarrow$ $C_{32} (3,4,13)$ and $C_{32} (4,5,11)$ are Type-1 isomorphic. 
\item [\rm (i2)]   $C_{32} (7(3,4,8,13))$ = $C_{32} (4,5,8,11)$.  $\Rightarrow$ $C_{32} (3,4,8,13)$ and $C_{32} (4,5,8,11)$ are Type-1 isomorphic.
\item [\rm (i3)]   $C_{32} (7(3,4,12,13))$ = $C_{32} (4,5,11,12)$.  
\\
$\Rightarrow$ $C_{32} (3,4,12,13)$ and $C_{32} (4,5,11,12)$ are Type-1 isomorphic.
\item [\rm (i4)]   $C_{32} (7(3,4,13,16))$ = $C_{32} (4,5,11,16)$.  
\\
$\Rightarrow$ $C_{32} (3,4,13,16)$ and $C_{32} (4,5,11,16)$ are Type-1 isomorphic. 
\item [\rm (i5)]   $C_{32} (7(3,8,12,13))$ = $C_{32} (5,8,11,12)$.  
\\
$\Rightarrow$ $C_{32} (3,8,12,13)$ and $C_{32} (5,8,11,12)$ are Type-1 isomorphic.
\item [\rm (i6)]   $C_{32} (7(3,12,13,16))$ = $C_{32} (5,11,12,16)$.  
\\
$\Rightarrow$ $C_{32} (3,12,13,16)$ and $C_{32} (5,11,12,16)$ are Type-1 isomorphic.
\item [\rm (i7)]   $C_{32} (7(3,4,8,12,13))$ = $C_{32} (4,5,8,11,12)$.  
\\
$\Rightarrow$ $C_{32} (3,4,8,12,13)$ and $C_{32} (4,5,8,11,12)$ are Type-1 isomorphic.
\item [\rm (i8)]   $C_{32} (7(3,4,8,13,16))$ = $C_{32} (4,5,8,11,16)$.  
\\
$\Rightarrow$ $C_{32} (3,4,8,13,16)$ and $C_{32} (4,5,8,11,16)$ are Type-1 isomorphic.
\item [\rm (i9)]   $C_{32} (7(3,4,12,13,16))$ = $C_{32} (4,5,11,12,16)$.  
\\
$\Rightarrow$ $C_{32} (3,4,12,13,16)$ and $C_{32} (4,5,11,12,16)$ are Type-1 isomorphic.
\item [\rm (i10)]  $C_{32} (7(3,8,12,13,16))$ = $C_{32} (5,8,11,12,16)$.  
\\
$\Rightarrow$ $C_{32} (3,8,12,13,16)$ and $C_{32} (5,8,11,12,16)$ are Type-1 isomorphic.
\item [\rm (i11)] $C_{32} (7(3,4,8,12,13,16))$ = $C_{32} (4,5,8,11,12,16)$.  
\\
$\Rightarrow$ $C_{32} (3,4,8,12,13,16)$ and $C_{32} (4,5,8,11,12,16)$ are Type-1 isomorphic. 
\item [\rm (j1)]   $C_{32} (7(3,8,13))$ = $C_{32} (5,8,11)$.  $\Rightarrow$ $C_{32} (3,8,13)$ and $C_{32} (5,8,11)$ are Type-1 isomorphic.
\item [\rm (j2)]   $C_{32} (7(3,8,13,16))$ = $C_{32} (5,8,11,16)$.  
\\
$\Rightarrow$ $C_{32} (3,8,13,16)$ and $C_{32} (5,8,11,16)$ are Type-1 isomorphic.
\item [\rm (k1)]   $C_{32} (7(3,12,13))$ = $C_{32} (5,11,12)$.  $\Rightarrow$ $C_{32} (3,12,13)$ and $C_{32} (5,11,12)$ are Type-1 isomorphic.
\item [\rm (l1)]   $C_{32} (7(3,13,16))$ = $C_{32} (5,11,16)$.  $\Rightarrow$ $C_{32} (3,13,16)$ and $C_{32} (5,11,16)$ are Type-1 isomorphic.
\end{enumerate}

\vspace{.2cm}
In the next problem, we find all Type-2 isomorphic circulant graphs of order 32 using remark \ref{r5} on the 8 pairs of Type-2 isomorphic circulant graphs, 

 (1) $C_{32} (1,2,15)$ and $C_{32} (2,7,9)$; (2) $C_{32}(1,6,15)$ and $C_{32}(6,7,9)$; 

(3) $C_{32}(1,10,15)$ and $C_{32}(7,9,10)$; (4) $C_{32}(1,14,15)$ and $C_{32}(7,9,14)$; 

(5) $C_{32}(2,3,13)$ and $C_{32}(2,5,11)$; (6)  $C_{32}(3,6,13)$ and $C_{32}(5,6,11)$; 

(7) $C_{32}(3,10,13)$ and $C_{32}(5,10,11)$; and (8) $C_{32}(3,13,14)$ and $C_{32}(5,11,14)$. 
\\
Note that the jump sizes of these 8 pairs of circulant graphs cover all jump sizes for all possible pairs of  Type-2 isomorphic circulant graphs $C_{32}(R)$ and $C_{32}(S)$ with $|R|$ = 3 = $|S|$.

\begin{prm}\quad \label{p3} {\rm Show that each pair of circulant graphs of order 32 given below are isomorphic; classify their type of isomorphism; and the number of pairs of Type-2 isomorphic circulant graphs is 384.
\begin{enumerate}
\item [\rm (a)] Pairs of circulant graphs of order 32 containing jump sizes 1, 15, 7 and 9.
\begin{enumerate}
\item [\rm (1)]  $C_{32} (1,2,15)$, $C_{32} (2,7,9)$;  

\item [\rm (2)] $C_{32} (1,4,15)$, $C_{32} (4,7,9)$; 

\item [\rm (3)] $C_{32} (1,6,15)$, $C_{32} (6,7,9)$; 
	
\item [\rm (4)] $C_{32} (1,8,15)$, $C_{32} (7,8,9)$; 

\item [\rm (5)] $C_{32} (1,10,15)$, $C_{32} (7,9,10)$;  
	
\item [\rm (6)] $C_{32} (1,12,15)$, $C_{32} (7,9,12)$;  
	
\item [\rm (7)] $C_{32} (1,14,15)$, $C_{32} (7,9,14)$; 
	
\item [\rm (8)] $C_{32} (1,15,16)$, $C_{32} (7,9,16)$;  

\item [\rm (9)] $C_{32} (1,2,4,15)$, $C_{32} (2,4,7,9)$;  
	
\item [\rm (10)] $C_{32} (1,2,6,15)$, $C_{32} (2,6,7,9)$;
	
\item [\rm (11)] $C_{32} (1,2,8,15)$, $C_{32} (2,7,8,9)$; 
	
\item [\rm (12)] $C_{32} (1,2,10,15)$, $C_{32} (2,7,9,10)$; 
	
\item [\rm (13)] $C_{32} (1,2,12,15)$, $C_{32} (2,7,9,12)$; 
	
\item [\rm (14)] $C_{32} (1,2,14,15)$, $C_{32} (2,7,9,14)$; 
	
\item [\rm (15)] $C_{32} (1,2,15,16)$, $C_{32} (2,7,9,16)$; 

\item [\rm (16)] $C_{32} (1,4,6,15)$, $C_{32} (4,6,7,9)$; 

\item [\rm (17)] $C_{32} (1,4,8,15)$, $C_{32} (4,7,8,9)$; 

\item [\rm (18)] $C_{32} (1,4,10,15)$, $C_{32} (4,7,9,10)$; 

\item [\rm (19)] $C_{32} (1,4,12,15)$, $C_{32} (4,7,9,12)$; 

\item [\rm (20)] $C_{32} (1,4,14,15)$, $C_{32} (4,7,9,14)$; 

\item [\rm (21)] $C_{32} (1,4,15,16)$, $C_{32} (4,7,9,16)$; 

\item [\rm (22)] $C_{32} (1,6,8,15)$, $C_{32} (6,7,8,9)$; 
	
\item [\rm (23)] $C_{32} (1,6,10,15)$, $C_{32} (6,7,9,10)$; 
	
\item [\rm (24)] $C_{32} (1,6,12,15)$, $C_{32} (6,7,9,12)$; 
	
\item [\rm (25)] $C_{32} (1,6,14,15)$, $C_{32} (6,7,9,14)$; 
	
\item [\rm (26)] $C_{32} (1,6,15,16)$, $C_{32} (6,7,9,16)$; 
	
\item [\rm (27)] $C_{32} (1,8,10,15)$, $C_{32} (7,8,9,10)$; 

\item [\rm (28)] $C_{32} (1,8,12,15)$, $C_{32} (7,8,9,12)$; 

\item [\rm (29)] $C_{32} (1,8,14,15)$, $C_{32} (7,8,9,14)$; 

\item [\rm (30)] $C_{32} (1,8,15,16)$, $C_{32} (7,8,9,16)$; 

\item [\rm (31)] $C_{32} (1,10,12,15)$, $C_{32} (7,9,10,12)$; 

\item [\rm (32)] $C_{32} (1,10,14,15)$, $C_{32} (7,9,10,14)$; 

\item [\rm (33)] $C_{32} (1,10,15,16)$, $C_{32} (7,9,10,16)$; 

\item [\rm (34)] $C_{32} (1,12,14,15)$, $C_{32} (7,9,12,14)$; 

\item [\rm (35)] $C_{32} (1,12,15,16)$, $C_{32} (7,9,12,16)$; 

\item [\rm (36)] $C_{32} (1,14,15,16)$, $C_{32} (7,9,14,16)$; 

\item [\rm (37)] $C_{32} (1,2,4,6,15)$, $C_{32} (2,4,6,7,9)$;

\item [\rm (38)] $C_{32} (1,2,4,8,15)$, $C_{32} (2,4,7,8,9)$;

\item [\rm (39)] $C_{32} (1,2,4,10,15)$, $C_{32} (2,4,7,9,10)$;
	
\item [\rm (40)] $C_{32} (1,2,4,12,15)$, $C_{32} (2,4,7,9,12)$; 

\item [\rm (41)] $C_{32} (1,2,4,14,15)$, $C_{32} (2,4,7,9,14)$;

\item [\rm (42)] $C_{32} (1,2,4,15,16)$, $C_{32} (2,4,7,9,16)$; 

\item [\rm (43)] $C_{32} (1,2,6,8,15)$, $C_{32} (2,6,7,8,9)$;

\item [\rm (44)] $C_{32} (1,2,6,10,15)$, $C_{32} (2,6,7,9,10)$;

\item [\rm (45)] $C_{32} (1,2,6,12,15)$, $C_{32} (2,6,7,9,12)$;

\item [\rm (46)] $C_{32} (1,2,6,14,15)$, $C_{32} (2,6,7,9,14)$;

\item [\rm (47)] $C_{32} (1,2,6,15,16)$, $C_{32} (2,6,7,9,16)$;

\item [\rm (48)] $C_{32} (1,2,8,10,15)$, $C_{32} (2,7,8,9,10)$;  

\item [\rm (49)] $C_{32} (1,2,8,12,15)$, $C_{32} (2,7,8,9,12)$;  

\item [\rm (50)] $C_{32} (1,2,8,14,15)$, $C_{32} (2,7,8,9,14)$; 

\item [\rm (51)] $C_{32} (1,2,8,15,16)$, $C_{32} (2,7,8,9,16)$;

\item [\rm (52)] $C_{32} (1,2,10,12,15)$, $C_{32} (2,7,9,10,12)$;  

\item [\rm (53)] $C_{32} (1,2,10,14,15)$, $C_{32} (2,7,9,10,14)$;  

\item [\rm (54)] $C_{32} (1,2,10,15,16)$, $C_{32} (2,7,9,10,16)$; 

\item [\rm (55)] $C_{32} (1,2,12,14,15)$, $C_{32} (2,7,9,12,14)$; 

\item [\rm (56)] $C_{32} (1,2,12,15,16)$, $C_{32} (2,7,9,12,16)$;

\item [\rm (57)] $C_{32} (1,2,14,15,16)$, $C_{32} (2,7,9,14,16)$; 

\item [\rm (58)] $C_{32} (1,4,6,8,15)$, $C_{32} (4,6,7,8,9)$; 

\item [\rm (59)] $C_{32} (1,4,6,10,15)$, $C_{32} (4,6,7,9,10)$;

\item [\rm (60)] $C_{32} (1,4,6,12,15)$, $C_{32} (4,6,7,9,12)$;

\item [\rm (61)] $C_{32} (1,4,6,14,15)$, $C_{32} (4,6,7,9,14)$;

\item [\rm (62)] $C_{32} (1,4,6,15,16)$, $C_{32} (4,6,7,9,16)$;

\item [\rm (63)] $C_{32} (1,4,8,10,15)$, $C_{32} (4,7,8,9,10)$;

\item [\rm (64)] $C_{32} (1,4,8,12,15)$, $C_{32} (4,7,8,9,12)$;

\item [\rm (65)] $C_{32} (1,4,8,14,15)$, $C_{32} (4,7,8,9,14)$;

\item [\rm (66)] $C_{32} (1,4,8,15,16)$, $C_{32} (4,7,8,9,16)$;

\item [\rm (67)] $C_{32} (1,4,10,12,15)$, $C_{32} (4,7,9,10,12)$;

\item [\rm (68)] $C_{32} (1,4,10,14,15)$, $C_{32} (4,7,9,10,14)$;

\item [\rm (69)] $C_{32} (1,4,10,15,16)$, $C_{32} (4,7,9,10,16)$;

\item [\rm (70)] $C_{32} (1,4,12,14,15)$, $C_{32} (4,7,9,12,14)$;

\item [\rm (71)] $C_{32} (1,4,12,15,16)$, $C_{32} (4,7,9,12,16)$;

\item [\rm (72)] $C_{32} (1,4,14,15,16)$, $C_{32} (4,7,9,14,16)$;

\item [\rm (73)] $C_{32} (1,6,8,10,15)$, $C_{32} (6,7,8,9,10)$;

\item [\rm (74)] $C_{32} (1,6,8,12,15)$, $C_{32} (6,7,8,9,12)$;

\item [\rm (75)] $C_{32} (1,6,8,14,15)$, $C_{32} (6,7,8,9,14)$; 

\item [\rm (76)] $C_{32} (1,6,8,15,16)$, $C_{32} (6,7,8,9,16)$; 

\item [\rm (77)] $C_{32} (1,6,10,12,15)$, $C_{32} (6,7,9,10,12)$;

\item [\rm (78)] $C_{32} (1,6,10,14,15)$, $C_{32} (6,7,9,10,14)$; 

\item [\rm (79)] $C_{32} (1,6,10,15,16)$, $C_{32} (6,7,9,10,16)$; 

\item [\rm (80)] $C_{32} (1,6,12,14,15)$, $C_{32} (6,7,9,12,14)$; 

\item [\rm (81)] $C_{32} (1,6,12,15,16)$, $C_{32} (6,7,9,12,16)$; 

\item [\rm (82)] $C_{32} (1,6,14,15,16)$, $C_{32} (6,7,9,14,16)$; 

\item [\rm (83)] $C_{32} (1,8,10,12,15)$, $C_{32} (7,8,9,10,12)$;

\item [\rm (84)] $C_{32} (1,8,10,14,15)$, $C_{32} (7,8,9,10,14)$; 

\item [\rm (85)] $C_{32} (1,8,10,15,16)$, $C_{32} (7,8,9,10,16)$;

\item [\rm (86)] $C_{32} (1,8,12,14,15)$, $C_{32} (7,8,9,12,14)$;

\item [\rm (87)] $C_{32} (1,8,12,15,16)$, $C_{32} (7,8,9,12,16)$;

\item [\rm (88)] $C_{32} (1,8,14,15,16)$, $C_{32} (7,8,9,14,16)$;

\item [\rm (89)] $C_{32} (1,10,12,14,15)$, $C_{32} (7,9,10,12,14)$; 

\item [\rm (90)] $C_{32} (1,10,12,15,16)$, $C_{32} (7,9,10,12,16)$; 

\item [\rm (91)] $C_{32} (1,10,14,15,16)$, $C_{32} (7,9,10,14,16)$; 

\item [\rm (92)] $C_{32} (1,12,14,15,16)$, $C_{32} (7,9,12,14,16)$; 

\item [\rm (93)] $C_{32} (1,2,4,6,8,15)$, $C_{32} (2,4,6,7,8,9)$; 

\item [\rm (94)] $C_{32} (1,2,4,6,10,15)$, $C_{32} (2,4,6,7,9,10)$; 

\item [\rm (95)] $C_{32} (1,2,4,6,12,15)$, $C_{32} (2,4,6,7,9,12)$; 

\item [\rm (96)] $C_{32} (1,2,4,6,14,15)$, $C_{32} (2,4,6,7,9,14)$; 

\item [\rm (97)] $C_{32} (1,2,4,6,15,16)$, $C_{32} (2,4,6,7,9,16)$; 

\item [\rm (98)] $C_{32} (1,2,4,8,10,15)$, $C_{32} (2,4,8,7,9,10)$; 

\item [\rm (99)] $C_{32} (1,2,4,8,12,15)$, $C_{32} (2,4,8,7,9,12)$; 

\item [\rm (100)] $C_{32} (1,2,4,8,14,15)$, $C_{32} (2,4,8,7,9,14)$; 

\item [\rm (101)] $C_{32} (1,2,4,8,15,16)$, $C_{32} (2,4,8,7,9,16)$; 

 \item [\rm (102)] $C_{32} (1,2,4,10,12,15)$, $C_{32} (2,4,7,9,10,12)$;

 \item [\rm (103)] $C_{32} (1,2,4,10,14,15)$, $C_{32} (2,4,7,9,10,14)$; 
 
 \item [\rm (104)] $C_{32} (1,2,4,10,15,16)$, $C_{32} (2,4,7,9,10,16)$; 

 \item [\rm (105)] $C_{32} (1,2,4,12,14,15)$, $C_{32} (2,4,7,9,12,14)$; 

\item [\rm (106)] $C_{32} (1,2,4,12,15,16)$, $C_{32} (2,4,7,9,12,16)$; 

\item [\rm (107)] $C_{32} (1,2,4,14,15,16)$, $C_{32} (2,4,7,9,14,16)$; 

\item [\rm (108)] $C_{32} (1,2,6,8,10,15)$, $C_{32} (2,6,7,8,9,10)$;

\item [\rm (109)] $C_{32} (1,2,6,8,12,15)$, $C_{32} (2,6,7,8,9,12)$; 

\item [\rm (110)] $C_{32} (1,2,6,8,14,15)$, $C_{32} (2,6,7,8,9,14)$; 

\item [\rm (111)] $C_{32} (1,2,6,8,15,16)$, $C_{32} (2,6,7,8,9,16)$; 

\item [\rm (112)] $C_{32} (1,2,6,10,12,15)$, $C_{32} (2,6,7,9,10,12)$; 

\item [\rm (113)] $C_{32} (1,2,6,10,14,15)$, $C_{32} (2,6,7,9,10,14)$; 

\item [\rm (114)] $C_{32} (1,2,6,10,15,16)$, $C_{32} (2,6,7,9,10,16)$; 

\item [\rm (115)] $C_{32} (1,2,6,12,14,15)$, $C_{32} (2,6,7,9,12,14)$; 

\item [\rm (116)] $C_{32} (1,2,6,12,15,16)$, $C_{32} (2,6,7,9,12,16)$;

\item [\rm (117)] $C_{32} (1,2,6,14,15,16)$, $C_{32} (2,6,7,9,14,16)$; 

\item [\rm (118)] $C_{32} (1,2,8,10,12,15)$, $C_{32} (2,7,8,9,10,12)$; 

\item [\rm (119)] $C_{32} (1,2,8,10,14,15)$, $C_{32} (2,7,8,9,10,14)$; 

\item [\rm (120)] $C_{32} (1,2,8,10,15,16)$, $C_{32} (2,7,8,9,10,16)$; 

\item [\rm (121)] $C_{32} (1,2,8,12,14,15)$, $C_{32} (2,7,8,9,12,14)$; 

\item [\rm (122)] $C_{32} (1,2,8,12,15,16)$, $C_{32} (2,7,8,9,12,16)$; 

\item [\rm (123)] $C_{32} (1,2,8,14,15,16)$, $C_{32} (2,7,8,9,14,16)$;

\item [\rm (124)] $C_{32} (1,2,10,12,14,15)$, $C_{32} (2,7,9,10,12,14)$; 

\item [\rm (125)] $C_{32} (1,2,10,12,15,16)$, $C_{32} (2,7,9,10,12,16)$; 

\item [\rm (126)] $C_{32} (1,2,10,14,15,16)$, $C_{32} (2,7,9,10,14,16)$; 
 
\item [\rm (127)] $C_{32} (1,2,12,14,15,16)$, $C_{32} (2,7,9,12,14,16)$; 

\item [\rm (128)] $C_{32} (1,4,6,8,10,15)$, $C_{32} (4,6,7,8,9,10)$; 

\item [\rm (129)] $C_{32} (1,4,6,8,12,15)$, $C_{32} (4,6,7,8,9,12)$; 

\item [\rm (130)] $C_{32} (1,4,6,8,14,15)$, $C_{32} (4,6,7,8,9,14)$; 

\item [\rm (131)] $C_{32} (1,4,6,8,15,16)$, $C_{32} (4,6,7,8,9,16)$; 

\item [\rm (132)] $C_{32} (1,4,6,10,12,15)$, $C_{32} (4,6,7,9,10,12)$;

\item [\rm (133)] $C_{32} (1,4,6,10,14,15)$, $C_{32} (4,6,7,9,10,14)$; 

 \item [\rm (134)] $C_{32} (1,4,6,10,15,16)$, $C_{32} (4,6,7,9,10,16)$; 

\item [\rm (135)] $C_{32} (1,4,6,12,14,15)$, $C_{32} (4,6,7,9,12,14)$; 

\item [\rm (136)] $C_{32} (1,4,6,12,15,16)$, $C_{32} (4,6,7,9,12,16)$; 

\item [\rm (137)] $C_{32} (1,4,6,14,15,16)$, $C_{32} (4,6,7,9,14,16)$; 

\item [\rm (138)] $C_{32} (1,4,8,10,12,15)$, $C_{32} (4,7,8,9,10,12)$; 

\item [\rm (139)] $C_{32} (1,4,8,10,14,15)$, $C_{32} (4,7,8,9,10,14)$; 

\item [\rm (140)] $C_{32} (1,4,8,10,15,16)$, $C_{32} (4,7,8,9,10,16)$; 

\item [\rm (141)] $C_{32} (1,4,8,12,14,15)$, $C_{32} (4,7,8,9,12,14)$; 

\item [\rm (142)] $C_{32} (1,4,8,12,15,16)$, $C_{32} (4,7,8,9,12,16)$; 

\item [\rm (143)] $C_{32} (1,4,8,14,15,16)$, $C_{32} (4,7,8,9,14,16)$; 

\item [\rm (144)] $C_{32} (1,4,10,12,14,15)$, $C_{32} (4,7,9,10,12,14)$; 

\item [\rm (145)] $C_{32} (1,4,10,12,15,16)$, $C_{32} (4,7,9,10,12,16)$; 

\item [\rm (146)] $C_{32} (1,4,10,14,15,16)$, $C_{32} (4,7,9,10,14,16)$; 

\item [\rm (147)] $C_{32} (1,4,12,14,15,16)$, $C_{32} (4,7,9,12,14,16)$; 

\item [\rm (148)] $C_{32} (1,6,8,10,12,15)$, $C_{32} (6,7,8,9,10,12)$; 

\item [\rm (149)] $C_{32} (1,6,8,10,14,15)$, $C_{32} (6,7,8,9,10,14)$; 

 \item [\rm (150)] $C_{32} (1,6,8,10,15,16)$, $C_{32} (6,7,8,9,10,16)$; 

\item [\rm (151)] $C_{32} (1,6,8,12,14,15)$, $C_{32} (6,7,8,9,12,14)$; 

 \item [\rm (152)] $C_{32} (1,6,8,12,15,16)$, $C_{32} (6,7,8,9,12,16)$; 

\item [\rm (153)] $C_{32} (1,6,8,14,15,16)$, $C_{32} (6,7,8,9,14,16)$; 

 \item [\rm (154)] $C_{32} (1,6,10,12,14,15)$, $C_{32} (6,7,9,10,12,14)$; 

 \item [\rm (155)] $C_{32} (1,6,10,12,15,16)$, $C_{32} (6,7,9,10,12,16)$; 

  \item [\rm (156)] $C_{32} (1,6,10,14,15,16)$, $C_{32} (6,7,9,10,14,16)$; 

 \item [\rm (157)] $C_{32} (1,6,12,14,15,16)$, $C_{32} (6,7,9,12,14,16)$; 

 \item [\rm (158)] $C_{32} (1,8,10,12,14,15)$, $C_{32} (7,8,9,10,12,14)$; 

 \item [\rm (159)] $C_{32} (1,8,10,12,15,16)$, $C_{32} (7,8,9,10,12,16)$; 

 \item [\rm (160)] $C_{32} (1,8,10,14,15,16)$, $C_{32} (7,8,9,10,14,16)$; 

 \item [\rm (161)] $C_{32} (1,8,12,14,15,16)$, $C_{32} (7,8,9,12,14,16)$; 

 \item [\rm (162)] $C_{32} (1,10,12,14,15,16)$, $C_{32} (7,9,10,12,14,16)$; 

 \item [\rm (163)] $C_{32} (1,2,4,6,8,10,15)$, $C_{32} (2,4,6,7,8,9,10)$; 

 \item [\rm (164)] $C_{32} (1,2,4,6,8,12,15)$, $C_{32} (2,4,6,7,8,9,12)$; 

 \item [\rm (165)] $C_{32} (1,2,4,6,8,14,15)$, $C_{32} (2,4,6,7,8,9,14)$; 

 \item [\rm (166)] $C_{32} (1,2,4,6,8,15,16)$, $C_{32} (2,4,6,7,8,9,16)$; 

  \item [\rm (167)] $C_{32} (1,2,4,6,10,12,15)$, $C_{32} (2,4,6,7,10,9,12)$; 

 \item [\rm (168)] $C_{32} (1,2,4,6,10,14,15)$, $C_{32} (2,4,6,7,9,10,14)$; 

 \item [\rm (169)] $C_{32} (1,2,4,6,10,15,16)$, $C_{32} (2,4,6,7,9,10,16)$; 

 \item [\rm (170)] $C_{32} (1,2,4,6,12,14,15)$, $C_{32} (2,4,6,7,9,12,14)$; 

 \item [\rm (171)] $C_{32} (1,2,4,6,12,15,16)$, $C_{32} (2,4,6,7,9,12,16)$; 

 \item [\rm (172)] $C_{32} (1,2,4,6,14,15,16)$, $C_{32} (2,4,6,7,9,14,16)$; 

 \item [\rm (173)] $C_{32} (1,2,4,8,10,12,15)$, $C_{32} (2,4,7,8,9,10,12)$; 

 \item [\rm (174)] $C_{32} (1,2,4,8,10,14,15)$, $C_{32} (2,4,7,8,9,10,14)$; 

 \item [\rm (175)] $C_{32} (1,2,4,8,10,15,16)$, $C_{32} (2,4,7,8,9,10,16)$; 

 \item [\rm (176)] $C_{32} (1,2,4,8,12,14,15)$, $C_{32} (2,4,7,8,9,12,14)$; 

 \item [\rm (177)] $C_{32} (1,2,4,8,12,15,16)$, $C_{32} (2,4,7,8,9,12,16)$; 

 \item [\rm (178)] $C_{32} (1,2,4,8,14,15,16)$, $C_{32} (2,4,7,8,9,14,16)$; 

 \item [\rm (179)] $C_{32} (1,2,4,10,12,14,15)$, $C_{32} (2,4,7,9,10,12,14)$; 

 \item [\rm (180)] $C_{32} (1,2,4,10,12,15,16)$, $C_{32} (2,4,7,9,10,12,16)$; 

 \item [\rm (181)] $C_{32} (1,2,4,10,14,15,16)$, $C_{32} (2,4,7,9,10,14,16)$; 

 \item [\rm (182)] $C_{32} (1,2,4,12,14,15,16)$, $C_{32} (2,4,7,9,12,14,16)$; 

 \item [\rm (183)] $C_{32} (1,2,6,8,10,12,15)$, $C_{32} (2,6,7,8,9,10,12)$; 

 \item [\rm (184)] $C_{32} (1,2,6,8,10,14,15)$, $C_{32} (2,6,7,8,9,10,14)$; 

 \item [\rm (185)] $C_{32} (1,2,6,8,10,15,16)$, $C_{32} (2,6,7,8,9,10,16)$; 

 \item [\rm (186)] $C_{32} (1,2,6,8,12,14,15)$, $C_{32} (2,6,7,8,9,12,14)$; 

\item [\rm (187)] $C_{32} (1,2,6,8,12,15,16)$, $C_{32} (2,6,7,8,9,12,16)$; 

\item [\rm (188)] $C_{32} (1,2,6,8,14,15,16)$, $C_{32} (2,6,7,8,9,14,16)$; 

\item [\rm (189)] $C_{32} (1,2,6,10,12,14,15)$, $C_{32} (2,6,7,9,10,12,14)$; 

\item [\rm (190)] $C_{32} (1,2,6,10,12,15,16)$, $C_{32} (2,6,7,9,10,12,16)$; 

\item [\rm (191)] $C_{32} (1,2,6,10,14,15,16)$, $C_{32} (2,6,7,9,10,14,16)$; 

\item [\rm (192)] $C_{32} (1,2,6,12,14,15,16)$, $C_{32} (2,6,7,9,12,14,16)$; 

\item [\rm (193)] $C_{32} (1,2,8,10,12,14,15)$, $C_{32} (2,7,8,9,10,12,14)$; 

\item [\rm (194)] $C_{32} (1,2,8,10,12,15,16)$, $C_{32} (2,7,8,9,10,12,16)$; 

\item [\rm (195)] $C_{32} (1,2,8,10,14,15,16)$, $C_{32} (2,7,8,9,10,14,16)$; 

\item [\rm (196)] $C_{32} (1,2,8,12,14,15,16)$, $C_{32} (2,7,8,9,12,14,16)$; 

\item [\rm (197)] $C_{32} (1,2,10,12,14,15,16)$, $C_{32} (2,7,9,10,12,14,16)$; 

\item [\rm (198)] $C_{32} (1,4,6,8,10,12,15)$, $C_{32} (4,6,7,8,9,10,12)$; 

\item [\rm (199)] $C_{32} (1,4,6,8,10,14,15)$, $C_{32} (4,6,7,8,9,10,14)$; 

\item [\rm (200)] $C_{32} (1,4,6,8,10,15,16)$, $C_{32} (4,6,7,8,9,10,16)$; 

\item [\rm (201)] $C_{32} (1,4,6,8,12,14,15)$, $C_{32} (4,6,7,8,9,12,14)$; 

\item [\rm (202)] $C_{32} (1,4,6,8,12,15,16)$, $C_{32} (4,6,7,8,9,12,16)$; 

\item [\rm (203)] $C_{32} (1,4,6,8,14,15,16)$, $C_{32} (4,6,7,8,9,14,16)$; 

\item [\rm (204)] $C_{32} (1,4,6,10,12,14,15)$, $C_{32} (4,6,7,9,10,12,14)$; 

\item [\rm (205)] $C_{32} (1,4,6,10,12,15,16)$, $C_{32} (4,6,7,9,10,12,16)$; 

\item [\rm (206)] $C_{32} (1,4,6,10,14,15,16)$, $C_{32} (4,6,7,9,10,14,16)$; 

\item [\rm (207)] $C_{32} (1,4,6,12,14,15,16)$, $C_{32} (4,6,7,9,12,14,16)$; 

\item [\rm (208)] $C_{32} (1,4,8,10,12,14,15)$, $C_{32} (4,7,8,9,10,12,14)$; 

\item [\rm (209)] $C_{32} (1,4,8,10,12,15,16)$, $C_{32} (4,7,8,9,10,12,16)$; 

\item [\rm (210)] $C_{32} (1,4,8,10,14,15,16)$, $C_{32} (4,7,8,9,10,14,16)$; 

\item [\rm (211)] $C_{32} (1,4,8,12,14,15,16)$, $C_{32} (4,7,8,9,12,14,16)$; 

\item [\rm (212)] $C_{32} (1,4,10,12,14,15,16)$, $C_{32} (4,7,9,10,12,14,16)$; 

\item [\rm (213)] $C_{32} (1,6,8,10,12,14,15)$, $C_{32} (6,7,8,9,10,12,14)$; 

\item [\rm (214)] $C_{32} (1,6,8,10,12,15,16)$, $C_{32} (6,7,8,9,10,12,16)$; 

\item [\rm (215)] $C_{32} (1,6,8,10,14,15,16)$, $C_{32} (6,7,8,9,10,14,16)$; 

\item [\rm (216)] $C_{32} (1,6,8,12,14,15,16)$, $C_{32} (6,7,8,9,12,14,16)$; 

\item [\rm (217)] $C_{32} (1,6,10,12,14,15,16)$, $C_{32} (6,7,9,10,12,14,16)$;

\item [\rm (218)] $C_{32} (1,8,10,12,14,15,16)$, $C_{32} (7,8,9,10,12,14,16)$; 
 
\item [\rm (219)] $C_{32} (1,2,4,6,8,10,12,15)$, $C_{32} (2,4,6,7,8,9,10,12)$; 

\item [\rm (220)] $C_{32} (1,2,4,6,8,10,14,15)$, $C_{32} (2,4,6,7,8,9,10,14)$; 

\item [\rm (221)] $C_{32} (1,2,4,6,8,10,15,16)$, $C_{32} (2,4,6,7,8,9,10,16)$;

\item [\rm (222)] $C_{32} (1,2,4,6,8,12,14,15)$, $C_{32} (2,4,6,7,8,9,12,14)$; 

\item [\rm (223)] $C_{32} (1,2,4,6,8,12,15,16)$, $C_{32} (2,4,6,7,8,9,12,16)$; 

\item [\rm (224)] $C_{32} (1,2,4,6,8,14,15,16)$, $C_{32} (2,4,6,7,8,9,14,16)$; 

 \item [\rm (225)] $C_{32} (1,2,4,6,10,12,14,15)$, $C_{32} (2,4,6,7,9,10,12,14)$; 

 \item [\rm (226)] $C_{32} (1,2,4,6,10,12,15,16)$, $C_{32} (2,4,6,7,9,10,12,16)$; 

 \item [\rm (227)] $C_{32} (1,2,4,6,10,14,15,16)$, $C_{32} (2,4,6,7,9,10,14,16)$; 

 \item [\rm (228)] $C_{32} (1,2,4,6,12,14,15,16)$, $C_{32} (2,4,6,7,9,12,14,16)$; 

 \item [\rm (229)] $C_{32} (1,2,4,8,10,12,14,15)$, $C_{32} (2,4,7,8,9,10,12,14)$; 

\item [\rm (230)] $C_{32} (1,2,4,8,10,12,15,16)$, $C_{32} (2,4,7,8,9,10,12,16)$; 

\item [\rm (231)] $C_{32} (1,2,4,8,10,14,15,16)$, $C_{32} (2,4,7,8,9,10,14,16)$; 

\item [\rm (232)] $C_{32} (1,2,4,8,12,14,15,16)$, $C_{32} (2,4,7,8,9,12,14,16)$; 

\item [\rm (233)] $C_{32} (1,2,4,10,12,14,15,16)$, $C_{32} (2,4,7,9,10,12,14,16)$; 

\item [\rm (234)] $C_{32} (1,2,6,8,10,12,14,15)$, $C_{32} (2,6,7,8,9,10,12,14)$; 

\item [\rm (235)] $C_{32} (1,2,6,8,10,12,15,16)$, $C_{32} (2,6,7,8,9,10,12,16)$; 

\item [\rm (236)] $C_{32} (1,2,6,8,10,14,15,16)$, $C_{32} (2,6,7,8,9,10,14,16)$; 

\item [\rm (237)] $C_{32} (1,2,6,8,12,14,15,16)$, $C_{32} (2,6,7,8,9,12,14,16)$; 

\item [\rm (238)] $C_{32} (1,2,6,10,12,14,15,16)$, $C_{32} (2,6,7,9,10,12,14,16)$; 

\item [\rm (239)] $C_{32} (1,2,8,10,12,14,15,16)$, $C_{32} (2,7,8,9,10,12,14,16)$; 

\item [\rm (240)] $C_{32} (1,4,6,8,10,12,14,15)$, $C_{32} (4,6,7,8,9,10,12,14)$; 

\item [\rm (241)] $C_{32} (1,4,6,8,10,12,15,16)$, $C_{32} (4,6,7,8,9,10,12,16)$; 

\item [\rm (242)] $C_{32} (1,4,6,8,10,14,15,16)$, $C_{32} (4,6,7,8,9,10,14,16)$; 

\item [\rm (243)] $C_{32} (1,4,6,8,12,14,15,16)$, $C_{32} (4,6,7,8,9,12,14,16)$; 

\item [\rm (244)] $C_{32} (1,4,6,10,12,14,15,16)$, $C_{32} (4,6,7,9,10,12,14,16)$; 

\item [\rm (245)] $C_{32} (1,4,8,10,12,14,15,16)$, $C_{32} (4,7,8,9,10,12,14,16)$; 

\item [\rm (246)] $C_{32} (1,6,8,10,12,14,15,16)$, $C_{32} (6,7,8,9,10,12,14,16)$; 

\item [\rm (247)] $C_{32} (1,2,4,6,8,10,12,14,15)$, $C_{32} (2,4,6,7,8,9,10,12,14)$; 

\item [\rm (248)] $C_{32} (1,2,4,6,8,10,12,15,16)$, $C_{32} (2,4,6,7,8,9,10,12,16)$; 

\item [\rm (249)] $C_{32} (1,2,4,6,8,10,14,15,16)$, $C_{32} (2,4,6,7,8,9,10,14,16)$; 

\item [\rm (250)] $C_{32} (1,2,4,6,8,12,14,15,16)$, $C_{32} (2,4,6,7,8,9,12,14,16)$; 

\item [\rm (251)] $C_{32} (1,2,4,6,10,12,14,15,16)$, $C_{32} (2,4,6,7,9,10,12,14,16)$; 

\item [\rm (252)] $C_{32} (1,2,4,8,10,12,14,15,16)$, $C_{32} (2,4,7,8,9,10,12,14,16)$; 

\item [\rm (253)] $C_{32} (1,2,6,8,10,12,14,15,16)$, $C_{32} (2,6,7,8,9,10,12,14,16)$; 

\item [\rm (254)] $C_{32} (1,4,6,8,10,12,14,15,16)$, $C_{32} (4,6,7,8,9,10,12,14,16)$; 

\item [\rm (255)] $C_{32} (1,2,4,6,8,10,12,14,15,16)$, $C_{32} (2,4,6,7,8,9,10,12,14,16)$. 
\end{enumerate} 

\vspace{.2cm}
\item [\rm (b)] Pairs of circulant graphs of order 32 containing jump sizes 3, 13, 5 and 11.
\begin{enumerate}
\item [\rm (1)]  $C_{32} (2,3,13)$, $C_{32} (2,5,11)$;  

\item [\rm (2)] $C_{32} (3,4,13)$, $C_{32} (4,5,11)$; 

\item [\rm (3)] $C_{32} (3,6,13)$, $C_{32} (5,6,11)$; 
	
\item [\rm (4)] $C_{32} (3,8,13)$, $C_{32} (5,8,11)$; 

\item [\rm (5)] $C_{32} (3,10,13)$, $C_{32} (5,10,11)$;  
	
\item [\rm (6)] $C_{32} (3,12,13)$, $C_{32} (5,11,12)$;  
	
\item [\rm (7)] $C_{32} (3,13,14)$, $C_{32} (5,11,14)$; 
	
\item [\rm (8)] $C_{32} (3,13,16)$, $C_{32} (5,11,16)$;  

\item [\rm (9)] $C_{32} (2,3,4,13)$, $C_{32} (2,4,5,11)$;  
	
\item [\rm (10)] $C_{32} (2,3,6,13)$, $C_{32} (2,5,6,11)$;
	
\item [\rm (11)] $C_{32} (2,3,8,13)$, $C_{32} (2,5,8,11)$; 
	
\item [\rm (12)] $C_{32} (2,3,10,13)$, $C_{32} (2,5,10,11)$; 
	
\item [\rm (13)] $C_{32} (2,3,12,13)$, $C_{32} (2,5,11,12)$; 
	
\item [\rm (14)] $C_{32} (2,3,13,14)$, $C_{32} (2,5,11,14)$; 
	
\item [\rm (15)] $C_{32} (2,3,13,16)$, $C_{32} (2,5,11,16)$; 

\item [\rm (16)] $C_{32} (3,4,6,13)$, $C_{32} (4,5,6,11)$; 

\item [\rm (17)] $C_{32} (3,4,8,13)$, $C_{32} (4,5,8,11)$; 

\item [\rm (18)] $C_{32} (3,4,10,13)$, $C_{32} (4,5,10,11)$; 

\item [\rm (19)] $C_{32} (3,4,12,13)$, $C_{32} (4,5,11,12)$; 

\item [\rm (20)] $C_{32} (3,4,13,14)$, $C_{32} (4,5,11,14)$; 

\item [\rm (21)] $C_{32} (3,4,13,16)$, $C_{32} (4,5,11,16)$; 

\item [\rm (22)] $C_{32} (3,6,8,13)$, $C_{32} (5,6,8,11)$; 
	
\item [\rm (23)] $C_{32} (3,6,10,13)$, $C_{32} (5,6,10,11)$; 
	
\item [\rm (24)] $C_{32} (3,6,12,13)$, $C_{32} (5,6,11,12)$; 
	
\item [\rm (25)] $C_{32} (3,6,13,14)$, $C_{32} (5,6,11,14)$; 
	
\item [\rm (26)] $C_{32} (3,6,13,16)$, $C_{32} (5,6,11,16)$; 
	
\item [\rm (27)] $C_{32} (3,8,10,13)$, $C_{32} (5,8,10,11)$; 

\item [\rm (28)] $C_{32} (3,8,12,13)$, $C_{32} (5,8,11,12)$; 

\item [\rm (29)] $C_{32} (3,8,13,14)$, $C_{32} (5,8,11,14)$; 

\item [\rm (30)] $C_{32} (3,8,13,16)$, $C_{32} (5,8,11,16)$; 

\item [\rm (31)] $C_{32} (3,10,12,13)$, $C_{32} (5,10,11,12)$; 

\item [\rm (32)] $C_{32} (3,10,13,14)$, $C_{32} (5,10,11,14)$; 

\item [\rm (33)] $C_{32} (3,10,13,16)$, $C_{32} (5,10,11,16)$; 

\item [\rm (34)] $C_{32} (3,12,13,14)$, $C_{32} (5,11,12,14)$; 

\item [\rm (35)] $C_{32} (3,12,13,16)$, $C_{32} (5,11,12,16)$; 

\item [\rm (36)] $C_{32} (3,13,14,16)$, $C_{32} (5,11,14,16)$; 

\item [\rm (37)] $C_{32} (2,3,4,6,13)$, $C_{32} (2,4,5,6,11)$;

\item [\rm (38)] $C_{32} (2,3,4,8,13)$, $C_{32} (2,4,5,8,11)$;

\item [\rm (39)] $C_{32} (2,3,4,10,13)$, $C_{32} (2,4,5,10,11)$;
	
\item [\rm (40)] $C_{32} (2,3,4,12,13)$, $C_{32} (2,4,5,11,12)$; 

\item [\rm (41)] $C_{32} (2,3,4,13,14)$, $C_{32} (2,4,5,11,14)$;

\item [\rm (42)] $C_{32} (2,3,4,13,16)$, $C_{32} (2,4,5,11,16)$; 

\item [\rm (43)] $C_{32} (2,3,6,8,13)$, $C_{32} (2,5,6,8,11)$;

\item [\rm (44)] $C_{32} (2,3,6,10,13)$, $C_{32} (2,5,6,10,11)$;

\item [\rm (45)] $C_{32} (2,3,6,12,13)$, $C_{32} (2,5,6,11,12)$;

\item [\rm (46)] $C_{32} (2,3,6,13,14)$, $C_{32} (2,5,6,11,14)$;

\item [\rm (47)] $C_{32} (2,3,6,13,16)$, $C_{32} (2,5,6,11,16)$;

\item [\rm (48)] $C_{32} (2,3,8,10,13)$, $C_{32} (2,5,8,10,11)$;  

\item [\rm (49)] $C_{32} (2,3,8,12,13)$, $C_{32} (2,5,8,11,12)$;  

\item [\rm (50)] $C_{32} (2,3,8,13,14)$, $C_{32} (2,5,8,11,14)$; 

\item [\rm (51)] $C_{32} (2,3,8,13,16)$, $C_{32} (2,5,8,11,16)$;

\item [\rm (52)] $C_{32} (2,3,10,12,13)$, $C_{32} (2,5,10,11,12)$;  

\item [\rm (53)] $C_{32} (2,3,10,13,14)$, $C_{32} (2,5,10,11,14)$;  

\item [\rm (54)] $C_{32} (2,3,10,13,16)$, $C_{32} (2,5,10,11,16)$; 

\item [\rm (55)] $C_{32} (2,3,12,13,14)$, $C_{32} (2,5,11,12,14)$; 

\item [\rm (56)] $C_{32} (2,3,12,13,16)$, $C_{32} (2,5,11,12,16)$;

\item [\rm (57)] $C_{32} (2,3,13,14,16)$, $C_{32} (2,5,11,14,16)$; 

\item [\rm (58)] $C_{32} (3,4,6,8,13)$, $C_{32} (4,5,6,8,11)$; 

\item [\rm (59)] $C_{32} (3,4,6,10,13)$, $C_{32} (4,5,6,10,11)$;

\item [\rm (60)] $C_{32} (3,4,6,12,13)$, $C_{32} (4,5,6,11,12)$;

\item [\rm (61)] $C_{32} (3,4,6,13,14)$, $C_{32} (4,5,6,11,14)$;

\item [\rm (62)] $C_{32} (3,4,6,13,16)$, $C_{32} (4,5,6,11,16)$;

\item [\rm (63)] $C_{32} (3,4,8,10,13)$, $C_{32} (4,5,8,10,11)$;

\item [\rm (64)] $C_{32} (3,4,8,12,13)$, $C_{32} (4,5,8,11,12)$;

\item [\rm (65)] $C_{32} (3,4,8,13,14)$, $C_{32} (4,5,8,11,14)$;

\item [\rm (66)] $C_{32} (3,4,8,13,16)$, $C_{32} (4,5,8,11,16)$;

\item [\rm (67)] $C_{32} (3,4,10,12,13)$, $C_{32} (4,5,10,11,12)$;

\item [\rm (68)] $C_{32} (3,4,10,13,14)$, $C_{32} (4,5,10,11,14)$;

\item [\rm (69)] $C_{32} (3,4,10,13,16)$, $C_{32} (4,5,10,11,16)$;

\item [\rm (70)] $C_{32} (3,4,12,13,14)$, $C_{32} (4,5,11,12,14)$;

\item [\rm (71)] $C_{32} (3,4,12,13,16)$, $C_{32} (4,5,11,12,16)$;

\item [\rm (72)] $C_{32} (3,4,13,14,16)$, $C_{32} (4,5,11,14,16)$;

\item [\rm (73)] $C_{32} (3,6,8,10,13)$, $C_{32} (5,6,8,10,11)$;

\item [\rm (74)] $C_{32} (3,6,8,12,13)$, $C_{32} (5,6,8,11,12)$;

\item [\rm (75)] $C_{32} (3,6,8,13,14)$, $C_{32} (5,6,8,11,14)$; 

\item [\rm (76)] $C_{32} (3,6,8,13,16)$, $C_{32} (5,6,8,11,16)$; 

\item [\rm (77)] $C_{32} (3,6,10,12,13)$, $C_{32} (5,6,10,11,12)$;

\item [\rm (78)] $C_{32} (3,6,10,13,14)$, $C_{32} (5,6,10,11,14)$; 

\item [\rm (79)] $C_{32} (3,6,10,13,16)$, $C_{32} (5,6,10,11,16)$; 

\item [\rm (80)] $C_{32} (3,6,12,13,14)$, $C_{32} (5,6,11,12,14)$; 

\item [\rm (81)] $C_{32} (3,6,12,13,16)$, $C_{32} (5,6,11,12,16)$; 

\item [\rm (82)] $C_{32} (3,6,13,14,16)$, $C_{32} (5,6,11,14,16)$; 

\item [\rm (83)] $C_{32} (3,8,10,12,13)$, $C_{32} (5,8,10,11,12)$;

\item [\rm (84)] $C_{32} (3,8,10,13,14)$, $C_{32} (5,8,10,11,14)$; 

\item [\rm (85)] $C_{32} (3,8,10,13,16)$, $C_{32} (5,8,10,11,16)$;

\item [\rm (86)] $C_{32} (3,8,12,13,14)$, $C_{32} (5,8,11,12,14)$;

\item [\rm (87)] $C_{32} (3,8,12,13,16)$, $C_{32} (5,8,11,12,16)$;

\item [\rm (88)] $C_{32} (3,8,13,14,16)$, $C_{32} (5,8,11,14,16)$;

\item [\rm (89)] $C_{32} (3,10,12,13,14)$, $C_{32} (5,10,11,12,14)$; 

\item [\rm (90)] $C_{32} (3,10,12,13,16)$, $C_{32} (5,10,11,12,16)$; 

\item [\rm (91)] $C_{32} (3,10,13,14,16)$, $C_{32} (5,10,11,14,16)$; 

\item [\rm (92)] $C_{32} (3,12,13,14,16)$, $C_{32} (5,11,12,14,16)$; 

\item [\rm (93)] $C_{32} (2,3,4,6,8,13)$, $C_{32} (2,4,5,6,8,11)$; 

\item [\rm (94)] $C_{32} (2,3,4,6,10,13)$, $C_{32} (2,4,5,6,10,11)$; 

\item [\rm (95)] $C_{32} (2,3,4,6,12,13)$, $C_{32} (2,4,5,6,11,12)$; 

\item [\rm (96)] $C_{32} (2,3,4,6,13,14)$, $C_{32} (2,4,5,6,11,14)$; 

\item [\rm (97)] $C_{32} (2,3,4,6,13,16)$, $C_{32} (2,4,5,6,11,16)$; 

\item [\rm (98)] $C_{32} (2,3,4,8,10,13)$, $C_{32} (2,4,8,5,10,11)$; 

\item [\rm (99)] $C_{32} (2,3,4,8,12,13)$, $C_{32} (2,4,8,5,11,12)$; 

\item [\rm (100)] $C_{32} (2,3,4,8,13,14)$, $C_{32} (2,4,8,5,11,14)$; 

\item [\rm (101)] $C_{32} (2,3,4,8,13,16)$, $C_{32} (2,4,8,5,11,16)$; 

 \item [\rm (102)] $C_{32} (2,3,4,10,12,13)$, $C_{32} (2,4,5,10,11,12)$;

 \item [\rm (103)] $C_{32} (2,3,4,10,13,14)$, $C_{32} (2,4,5,10,11,14)$; 
 
 \item [\rm (104)] $C_{32} (2,3,4,10,13,16)$, $C_{32} (2,4,5,10,11,16)$; 

 \item [\rm (105)] $C_{32} (2,3,4,12,13,14)$, $C_{32} (2,4,5,11,12,14)$; 

\item [\rm (106)] $C_{32} (2,3,4,12,13,16)$, $C_{32} (2,4,5,11,12,16)$; 

\item [\rm (107)] $C_{32} (2,3,4,13,14,16)$, $C_{32} (2,4,5,11,14,16)$; 

\item [\rm (108)] $C_{32} (2,3,6,8,10,13)$, $C_{32} (2,5,6,8,10,11)$;

\item [\rm (109)] $C_{32} (2,3,6,8,12,13)$, $C_{32} (2,5,6,8,11,12)$; 

\item [\rm (110)] $C_{32} (2,3,6,8,13,14)$, $C_{32} (2,5,6,8,11,14)$; 

\item [\rm (111)] $C_{32} (2,3,6,8,13,16)$, $C_{32} (2,5,6,8,11,16)$; 

\item [\rm (112)] $C_{32} (2,3,6,10,12,13)$, $C_{32} (2,5,6,10,11,12)$; 

\item [\rm (113)] $C_{32} (2,3,6,10,13,14)$, $C_{32} (2,5,6,10,11,14)$; 

\item [\rm (114)] $C_{32} (2,3,6,10,13,16)$, $C_{32} (2,5,6,10,11,16)$; 

\item [\rm (115)] $C_{32} (2,3,6,12,13,14)$, $C_{32} (2,5,6,11,12,14)$; 

\item [\rm (116)] $C_{32} (2,3,6,12,13,16)$, $C_{32} (2,5,6,11,12,16)$;

\item [\rm (117)] $C_{32} (2,3,6,13,14,16)$, $C_{32} (2,5,6,11,14,16)$; 

\item [\rm (118)] $C_{32} (2,3,8,10,12,13)$, $C_{32} (2,5,8,10,11,12)$; 

\item [\rm (119)] $C_{32} (2,3,8,10,13,14)$, $C_{32} (2,5,8,10,11,14)$; 

\item [\rm (120)] $C_{32} (2,3,8,10,13,16)$, $C_{32} (2,5,8,10,11,16)$; 

\item [\rm (121)] $C_{32} (2,3,8,12,13,14)$, $C_{32} (2,5,8,11,12,14)$; 

\item [\rm (122)] $C_{32} (2,3,8,12,13,16)$, $C_{32} (2,5,8,11,12,16)$; 

\item [\rm (123)] $C_{32} (2,3,8,13,14,16)$, $C_{32} (2,5,8,11,14,16)$;

\item [\rm (124)] $C_{32} (2,3,10,12,13,14)$, $C_{32} (2,5,10,11,12,14)$; 

\item [\rm (125)] $C_{32} (2,3,10,12,13,16)$, $C_{32} (2,5,10,11,12,16)$; 

\item [\rm (126)] $C_{32} (2,3,10,13,14,16)$, $C_{32} (2,5,10,11,14,16)$; 
 
\item [\rm (127)] $C_{32} (2,3,12,13,14,16)$, $C_{32} (2,5,11,12,14,16)$; 

\item [\rm (128)] $C_{32} (3,4,6,8,10,13)$, $C_{32} (4,5,6,8,10,11)$; 

\item [\rm (129)] $C_{32} (3,4,6,8,12,13)$, $C_{32} (4,5,6,8,11,12)$; 

\item [\rm (130)] $C_{32} (3,4,6,8,13,14)$, $C_{32} (4,5,6,8,11,14)$; 

\item [\rm (131)] $C_{32} (3,4,6,8,13,16)$, $C_{32} (4,5,6,8,11,16)$; 

\item [\rm (132)] $C_{32} (3,4,6,10,12,13)$, $C_{32} (4,5,6,10,11,12)$;

\item [\rm (133)] $C_{32} (3,4,6,10,13,14)$, $C_{32} (4,5,6,10,11,14)$; 

 \item [\rm (134)] $C_{32} (3,4,6,10,13,16)$, $C_{32} (4,5,6,10,11,16)$; 

\item [\rm (135)] $C_{32} (3,4,6,12,13,14)$, $C_{32} (4,5,6,11,12,14)$; 

\item [\rm (136)] $C_{32} (3,4,6,12,13,16)$, $C_{32} (4,5,6,11,12,16)$; 

\item [\rm (137)] $C_{32} (3,4,6,13,14,16)$, $C_{32} (4,5,6,11,14,16)$; 

\item [\rm (138)] $C_{32} (3,4,8,10,12,13)$, $C_{32} (4,5,8,10,11,12)$; 

\item [\rm (139)] $C_{32} (3,4,8,10,13,14)$, $C_{32} (4,5,8,10,11,14)$; 

\item [\rm (140)] $C_{32} (3,4,8,10,13,16)$, $C_{32} (4,5,8,10,11,16)$; 

\item [\rm (141)] $C_{32} (3,4,8,12,13,14)$, $C_{32} (4,5,8,11,12,14)$; 

\item [\rm (142)] $C_{32} (3,4,8,12,13,16)$, $C_{32} (4,5,8,11,12,16)$; 

\item [\rm (143)] $C_{32} (3,4,8,13,14,16)$, $C_{32} (4,5,8,11,14,16)$; 

\item [\rm (144)] $C_{32} (3,4,10,12,13,14)$, $C_{32} (4,5,10,11,12,14)$; 

\item [\rm (145)] $C_{32} (3,4,10,12,13,16)$, $C_{32} (4,5,10,11,12,16)$; 

\item [\rm (146)] $C_{32} (3,4,10,13,14,16)$, $C_{32} (4,5,10,11,14,16)$; 

\item [\rm (147)] $C_{32} (3,4,12,13,14,16)$, $C_{32} (4,5,11,12,14,16)$; 

\item [\rm (148)] $C_{32} (3,6,8,10,12,13)$, $C_{32} (5,6,8,10,11,12)$; 

\item [\rm (149)] $C_{32} (3,6,8,10,13,14)$, $C_{32} (5,6,8,10,11,14)$; 

 \item [\rm (150)] $C_{32} (3,6,8,10,13,16)$, $C_{32} (5,6,8,10,11,16)$; 

\item [\rm (151)] $C_{32} (3,6,8,12,13,14)$, $C_{32} (5,6,8,11,12,14)$; 

 \item [\rm (152)] $C_{32} (3,6,8,12,13,16)$, $C_{32} (5,6,8,11,12,16)$; 

\item [\rm (153)] $C_{32} (3,6,8,13,14,16)$, $C_{32} (5,6,8,11,14,16)$; 

 \item [\rm (154)] $C_{32} (3,6,10,12,13,14)$, $C_{32} (5,6,10,11,12,14)$; 

 \item [\rm (155)] $C_{32} (3,6,10,12,13,16)$, $C_{32} (5,6,10,11,12,16)$; 

  \item [\rm (156)] $C_{32} (3,6,10,13,14,16)$, $C_{32} (5,6,10,11,14,16)$; 

 \item [\rm (157)] $C_{32} (3,6,12,13,14,16)$, $C_{32} (5,6,11,12,14,16)$; 

 \item [\rm (158)] $C_{32} (3,8,10,12,13,14)$, $C_{32} (5,8,10,11,12,14)$; 

 \item [\rm (159)] $C_{32} (3,8,10,12,13,16)$, $C_{32} (5,8,10,11,12,16)$; 

 \item [\rm (160)] $C_{32} (3,8,10,13,14,16)$, $C_{32} (5,8,10,11,14,16)$; 

 \item [\rm (161)] $C_{32} (3,8,12,13,14,16)$, $C_{32} (5,8,11,12,14,16)$; 

 \item [\rm (162)] $C_{32} (3,10,12,13,14,16)$, $C_{32} (5,10,11,12,14,16)$; 

 \item [\rm (163)] $C_{32} (2,3,4,6,8,10,13)$, $C_{32} (2,4,5,6,8,10,11)$; 

 \item [\rm (164)] $C_{32} (2,3,4,6,8,12,13)$, $C_{32} (2,4,5,6,8,11,12)$; 

 \item [\rm (165)] $C_{32} (2,3,4,6,8,13,14)$, $C_{32} (2,4,5,6,8,11,14)$; 

 \item [\rm (166)] $C_{32} (2,3,4,6,8,13,16)$, $C_{32} (2,4,5,6,8,11,16)$; 

  \item [\rm (167)] $C_{32} (2,3,4,6,10,12,13)$, $C_{32} (2,4,5,6,10,11,12)$; 

 \item [\rm (168)] $C_{32} (2,3,4,6,10,13,14)$, $C_{32} (2,4,5,6,10,11,14)$; 

 \item [\rm (169)] $C_{32} (2,3,4,6,10,13,16)$, $C_{32} (2,4,5,6,10,11,16)$; 

 \item [\rm (170)] $C_{32} (2,3,4,6,12,13,14)$, $C_{32} (2,4,5,6,11,12,14)$; 

 \item [\rm (171)] $C_{32} (2,3,4,6,12,13,16)$, $C_{32} (2,4,5,6,11,12,16)$; 

 \item [\rm (172)] $C_{32} (2,3,4,6,13,14,16)$, $C_{32} (2,4,5,6,11,14,16)$; 

 \item [\rm (173)] $C_{32} (2,3,4,8,10,12,13)$, $C_{32} (2,4,5,8,10,11,12)$; 

 \item [\rm (174)] $C_{32} (2,3,4,8,10,13,14)$, $C_{32} (2,4,5,8,10,11,14)$; 

 \item [\rm (175)] $C_{32} (2,3,4,8,10,13,16)$, $C_{32} (2,4,5,8,10,11,16)$; 

 \item [\rm (176)] $C_{32} (2,3,4,8,12,13,14)$, $C_{32} (2,4,5,8,11,12,14)$; 

 \item [\rm (177)] $C_{32} (2,3,4,8,12,13,16)$, $C_{32} (2,4,5,8,11,12,16)$; 

 \item [\rm (178)] $C_{32} (2,3,4,8,13,14,16)$, $C_{32} (2,4,5,8,11,14,16)$; 

 \item [\rm (179)] $C_{32} (2,3,4,10,12,13,14)$, $C_{32} (2,4,5,10,11,12,14)$; 

 \item [\rm (180)] $C_{32} (2,3,4,10,12,13,16)$, $C_{32} (2,4,5,10,11,12,16)$; 

 \item [\rm (181)] $C_{32} (2,3,4,10,13,14,16)$, $C_{32} (2,4,5,10,11,14,16)$; 

 \item [\rm (182)] $C_{32} (2,3,4,12,13,14,16)$, $C_{32} (2,4,5,11,12,14,16)$; 

 \item [\rm (183)] $C_{32} (2,3,6,8,10,12,13)$, $C_{32} (2,5,6,8,10,11,12)$; 

 \item [\rm (184)] $C_{32} (2,3,6,8,10,13,14)$, $C_{32} (2,5,6,8,10,11,14)$; 

 \item [\rm (185)] $C_{32} (2,3,6,8,10,13,16)$, $C_{32} (2,5,6,8,10,11,16)$; 

 \item [\rm (186)] $C_{32} (2,3,6,8,12,13,14)$, $C_{32} (2,5,6,8,11,12,14)$; 

\item [\rm (187)] $C_{32} (2,3,6,8,12,13,16)$, $C_{32} (2,5,6,8,11,12,16)$; 

\item [\rm (188)] $C_{32} (2,3,6,8,13,14,16)$, $C_{32} (2,5,6,8,11,14,16)$; 

\item [\rm (189)] $C_{32} (2,3,6,10,12,13,14)$, $C_{32} (2,5,6,10,11,12,14)$; 

\item [\rm (190)] $C_{32} (2,3,6,10,12,13,16)$, $C_{32} (2,5,6,10,11,12,16)$; 

\item [\rm (191)] $C_{32} (2,3,6,10,13,14,16)$, $C_{32} (2,5,6,10,11,14,16)$; 

\item [\rm (192)] $C_{32} (2,3,6,12,13,14,16)$, $C_{32} (2,5,6,11,12,14,16)$; 

\item [\rm (193)] $C_{32} (2,3,8,10,12,13,14)$, $C_{32} (2,5,8,10,11,12,14)$; 

\item [\rm (194)] $C_{32} (2,3,8,10,12,13,16)$, $C_{32} (2,5,8,10,11,12,16)$; 

\item [\rm (195)] $C_{32} (2,3,8,10,13,14,16)$, $C_{32} (2,5,8,10,11,14,16)$; 

\item [\rm (196)] $C_{32} (2,3,8,12,13,14,16)$, $C_{32} (2,5,8,11,12,14,16)$; 

\item [\rm (197)] $C_{32} (2,3,10,12,13,14,16)$, $C_{32} (2,5,10,11,12,14,16)$; 

\item [\rm (198)] $C_{32} (3,4,6,8,10,12,13)$, $C_{32} (4,5,6,8,10,11,12)$; 

\item [\rm (199)] $C_{32} (3,4,6,8,10,13,14)$, $C_{32} (4,5,6,8,10,11,14)$; 

\item [\rm (200)] $C_{32} (3,4,6,8,10,13,16)$, $C_{32} (4,5,6,8,10,11,16)$; 

\item [\rm (201)] $C_{32} (3,4,6,8,12,13,14)$, $C_{32} (4,5,6,8,11,12,14)$; 

\item [\rm (202)] $C_{32} (3,4,6,8,12,13,16)$, $C_{32} (4,5,6,8,11,12,16)$; 

\item [\rm (203)] $C_{32} (3,4,6,8,13,14,16)$, $C_{32} (4,5,6,8,11,14,16)$; 

\item [\rm (204)] $C_{32} (3,4,6,10,12,13,14)$, $C_{32} (4,5,6,10,11,12,14)$; 

\item [\rm (205)] $C_{32} (3,4,6,10,12,13,16)$, $C_{32} (4,5,6,10,11,12,16)$; 

\item [\rm (206)] $C_{32} (3,4,6,10,13,14,16)$, $C_{32} (4,5,6,10,11,14,16)$; 

\item [\rm (207)] $C_{32} (3,4,6,12,13,14,16)$, $C_{32} (4,5,6,11,12,14,16)$; 

\item [\rm (208)] $C_{32} (3,4,8,10,12,13,14)$, $C_{32} (4,5,8,10,11,12,14)$; 

\item [\rm (209)] $C_{32} (3,4,8,10,12,13,16)$, $C_{32} (4,5,8,10,11,12,16)$; 

\item [\rm (210)] $C_{32} (3,4,8,10,13,14,16)$, $C_{32} (4,5,8,10,11,14,16)$; 

\item [\rm (211)] $C_{32} (3,4,8,12,13,14,16)$, $C_{32} (4,5,8,11,12,14,16)$; 

\item [\rm (212)] $C_{32} (3,4,10,12,13,14,16)$, $C_{32} (4,5,10,11,12,14,16)$; 

\item [\rm (213)] $C_{32} (3,6,8,10,12,13,14)$, $C_{32} (5,6,8,10,11,12,14)$; 

\item [\rm (214)] $C_{32} (3,6,8,10,12,13,16)$, $C_{32} (5,6,8,10,11,12,16)$; 

\item [\rm (215)] $C_{32} (3,6,8,10,13,14,16)$, $C_{32} (5,6,8,10,11,14,16)$; 

\item [\rm (216)] $C_{32} (3,6,8,12,13,14,16)$, $C_{32} (5,6,8,11,12,14,16)$; 

\item [\rm (217)] $C_{32} (3,6,10,12,13,14,16)$, $C_{32} (5,6,10,11,12,14,16)$;

\item [\rm (218)] $C_{32} (3,8,10,12,13,14,16)$, $C_{32} (7,8,10,11,12,14,16)$; 
 
\item [\rm (219)] $C_{32} (2,3,4,6,8,10,12,13)$, $C_{32} (2,4,5,6,8,10,11,12)$; 

\item [\rm (220)] $C_{32} (2,3,4,6,8,10,13,14)$, $C_{32} (2,4,5,6,8,10,11,14)$; 

\item [\rm (221)] $C_{32} (2,3,4,6,8,10,13,16)$, $C_{32} (2,4,5,6,8,10,11,16)$;

\item [\rm (222)] $C_{32} (2,3,4,6,8,12,13,14)$, $C_{32} (2,4,5,6,8,11,12,14)$; 

\item [\rm (223)] $C_{32} (2,3,4,6,8,12,13,16)$, $C_{32} (2,4,5,6,8,11,12,16)$; 

\item [\rm (224)] $C_{32} (2,3,4,6,8,13,14,16)$, $C_{32} (2,4,5,6,8,11,14,16)$; 

 \item [\rm (225)] $C_{32} (2,3,4,6,10,12,13,14)$, $C_{32} (2,4,5,6,10,11,12,14)$; 

 \item [\rm (226)] $C_{32} (2,3,4,6,10,12,13,16)$, $C_{32} (2,4,5,6,10,11,12,16)$; 

 \item [\rm (227)] $C_{32} (2,3,4,6,10,13,14,16)$, $C_{32} (2,4,5,6,10,11,14,16)$; 

 \item [\rm (228)] $C_{32} (2,3,4,6,12,13,14,16)$, $C_{32} (2,4,5,6,11,12,14,16)$; 

 \item [\rm (229)] $C_{32} (2,3,4,8,10,12,13,14)$, $C_{32} (2,4,5,8,10,11,12,14)$; 

\item [\rm (230)] $C_{32} (2,3,4,8,10,12,13,16)$, $C_{32} (2,4,5,8,10,11,12,16)$; 

\item [\rm (231)] $C_{32} (2,3,4,8,10,13,14,16)$, $C_{32} (2,4,5,8,10,11,14,16)$; 

\item [\rm (232)] $C_{32} (2,3,4,8,12,13,14,16)$, $C_{32} (2,4,5,8,11,12,14,16)$; 

\item [\rm (233)] $C_{32} (2,3,4,10,12,13,14,16)$, $C_{32} (2,4,5,10,11,12,14,16)$; 

\item [\rm (234)] $C_{32} (2,3,6,8,10,12,13,14)$, $C_{32} (2,5,6,8,10,11,12,14)$; 

\item [\rm (235)] $C_{32} (2,3,6,8,10,12,13,16)$, $C_{32} (2,5,6,8,10,11,12,16)$; 

\item [\rm (236)] $C_{32} (2,3,6,8,10,13,14,16)$, $C_{32} (2,5,6,8,10,11,14,16)$; 

\item [\rm (237)] $C_{32} (2,3,6,8,12,13,14,16)$, $C_{32} (2,5,6,8,11,12,14,16)$; 

\item [\rm (238)] $C_{32} (2,3,6,10,12,13,14,16)$, $C_{32} (2,5,6,10,11,12,14,16)$; 

\item [\rm (239)] $C_{32} (2,3,8,10,12,13,14,16)$, $C_{32} (2,5,8,10,11,12,14,16)$; 

\item [\rm (240)] $C_{32} (3,4,6,8,10,12,13,14)$, $C_{32} (4,5,6,8,10,11,12,14)$; 

\item [\rm (241)] $C_{32} (3,4,6,8,10,12,13,16)$, $C_{32} (4,5,6,8,10,11,12,16)$; 

\item [\rm (242)] $C_{32} (3,4,6,8,10,13,14,16)$, $C_{32} (4,5,6,8,10,11,14,16)$; 

\item [\rm (243)] $C_{32} (3,4,6,8,12,13,14,16)$, $C_{32} (4,5,6,8,11,12,14,16)$; 

\item [\rm (244)] $C_{32} (3,4,6,10,12,13,14,16)$, $C_{32} (4,5,6,10,11,12,14,16)$; 

\item [\rm (245)] $C_{32} (3,4,8,10,12,13,14,16)$, $C_{32} (4,5,8,10,11,12,14,16)$; 

\item [\rm (246)] $C_{32} (3,6,8,10,12,13,14,16)$, $C_{32} (5,6,8,10,11,12,14,16)$; 

\item [\rm (247)] $C_{32} (2,3,4,6,8,10,12,13,14)$, $C_{32} (2,4,5,6,8,10,11,12,14)$; 

\item [\rm (248)] $C_{32} (2,3,4,6,8,10,12,13,16)$, $C_{32} (2,4,5,6,8,10,11,12,16)$; 

\item [\rm (249)] $C_{32} (2,3,4,6,8,10,13,14,16)$, $C_{32} (2,4,5,6,8,10,11,14,16)$; 

\item [\rm (250)] $C_{32} (2,3,4,6,8,12,13,14,16)$, $C_{32} (2,4,5,6,8,11,12,14,16)$; 

\item [\rm (251)] $C_{32} (2,3,4,6,10,12,13,14,16)$, $C_{32} (2,4,5,6,10,11,12,14,16)$; 

\item [\rm (252)] $C_{32} (2,3,4,8,10,12,13,14,16)$, $C_{32} (2,4,5,8,10,11,12,14,16)$; 

\item [\rm (253)] $C_{32} (2,3,6,8,10,12,13,14,16)$, $C_{32} (2,5,6,8,10,11,12,14,16)$; 

\item [\rm (254)] $C_{32} (3,4,6,8,10,12,13,14,16)$, $C_{32} (4,5,6,8,10,11,12,14,16)$; 

\item [\rm (255)] $C_{32} (2,3,4,6,8,10,12,13,14,16)$, $C_{32} (2,4,5,6,8,10,11,12,14,16)$. 
\end{enumerate} 
\end{enumerate} } 
\end{prm}
\noindent
{\bf Solution.}\quad Here, we consider Type-1 and Type-2 isomorphisms of circulant graphs of the form $C_{32}(R)$ and so the possible values of $m > 1$ for the existence of Type-2 isomorphism w.r.t. $m$ of $C_{32}(R)$ $\ni$ $m$ is a divisor of $\gcd(n, r)$ = $\gcd(32, r)$ and $m^3$ divides $n$ = 32 is $m$ = 2. Also, we have $m$ = 2 = $\gcd(32, 2)$ = $\gcd(32, 6)$ = $\gcd(32, 10)$ = $\gcd(32, 14)$, 4 = $\gcd(32, 4)$ = $\gcd(32, 12)$, 8 = $\gcd(32, 8)$ and 16 = $\gcd(32, 16)$. 

At first, we consider pairs of circulant graphs of order 32 with jump sizes 1, 15, 7 and 9 and then we consider pairs of circulant graphs of order 32 with jump sizes 3, 13, 5 and 11. 

\vspace{.2cm}
\noindent
{\bf Case (a)}\quad  In problem \ref{p3.5}, we proved that circulant graphs 

\vspace{.1cm}
(1) $C_{32} (1, 2, 15)$ and $C_{32} (2, 7, 9)$ are isomorphic of Type-2 w.r.t. $m$ = 2. 
\\
Using remark \ref{r5} in this pair of Type-2 isomorphic circulant graphs, we obtain 234 pairs of isomorphic circulant graphs as given in the problem. Also, for a given circulant graph $C_n(R)$, if all $C_n(S)$ $\ni$ $C_n(S)$ = $\theta_{n,m,t}(C_{n}(R))$ for some $t$ and $C_n(S)\in T1_n(C_n(R))$, $1 \leq t \leq \frac{n}{m}-1$, then $C_n(R)$ has no isomorphic circulant graph of Type-2  w.r.t. $m$ where $r\in R,S$ and $m > 1$ and $m^3$ are divisors of $\gcd(n, r)$ and $n$, respectively. 

Solutions to all the cases in this problem are similar, in cases of proving Type-1 isomorphism and also in the cases of Type-2 isomorphism, and to simplify our work, we present important values related to all cases related to Type-2 isomorphism in Table 1 to Table 16. In each table, in the column of `T1 or T2' corresponds to whether the pair of isomorphic circulant graphs are `Type-1 isomorphic or Type-2 isomorphic'. Even though some pairs of circulant graphs are already covered in problem \ref{p2}, we consider all the pairs of isomorphic circulant graphs. 

\begin{table}
	\caption{ Finding $\theta_{32,2,4}(C_{32}(R))$ and $T1_{32}(C_{32}(R))$ as $\theta_{32,2,8}(C_{32}(R))$ = $C_{32}(R)$.}
\begin{center}
\scalebox{.8}{
\begin{tabular}{||c||c|c|c|c||}\hline \hline 
  & & & &  \\
S. No. & $C_{32}(R)$ &  $\theta_{32,2,4t}(C_{32}(R))$, $t$ = 1,2  & $C_{32}(S)\in T1_{24}(C_{32}(R))$ & T1 or T2  \\

 & & &   &  \\ \hline \hline

1.   & $C_{32}(1,2,15)$ &  = $C_{32}(2,7,9)$  & $C_{32}(1,2,15)$, $C_{32}(3,6,13)$,   &   \\ 

     & & &  $C_{32}(5,10,11)$, $C_{32}(7,9,14)$.  &  T2 \\

2.   & $C_{32}(1,4,15)$ &  = $C_{32}(4,7,9)$  & $C_{32}(1,4,15)$, $C_{32}(3,12,13)$,   &   \\ 

     & & &  $C_{32}(5,11,12)$, $C_{32}(4,7,9)$.  &  T1 \\

3.   & $C_{32}(1,6,15)$ &  = $C_{32}(6,7,9)$  & $C_{32}(1,6,15)$, $C_{32}(3,13,14)$,    &   \\
\hline \hline
\end{tabular}}
\end{center}
\end{table} 

\begin{table}
	\caption{ Finding $\theta_{32,2,4}(C_{32}(R))$ and $T1_{32}(C_{32}(R))$ as $\theta_{32,2,8}(C_{32}(R))$ = $C_{32}(R)$.}
\begin{center}
\scalebox{.8}{
 \begin{tabular}{||c||c|c|c|c||}\hline \hline 
  & & & &  \\
S. No. & $C_{32}(R)$ &  $\theta_{32,2,4t}(C_{32}(R))$, $t$ = 1,2  & $C_{32}(S)\in T1_{24}(C_{32}(R))$ & T1 or T2  \\

 & & &   &  \\ \hline \hline

31.  & $C_{32}(1,10,12,15)$ &  = $C_{32}(7,9,10,12)$  & $C_{32}(1,10,12,15)$, $C_{32}(2,3,4,13)$,   &   \\ 

     & & &  $C_{32}(4,5,11,14)$, $C_{32}(6,7,9,12)$.  & T2 \\      

32.  & $C_{32}(1,10,14,15)$ &  = $C_{32}(7,9,10,14)$  & $C_{32}(1,10,14,15)$, $C_{32}(2,3,10,13)$,   &   \\ 

     & & &  $C_{32}(5,6,11,14)$, $C_{32}(2,6,7,9)$.  & T2 \\   

33.  & $C_{32}(1,10,15,16)$ &  = $C_{32}(7,9,10,16)$  & $C_{32}(1,10,15,16)$, $C_{32}(2,3,13,16)$,   &   \\
\hline \hline
\end{tabular}}
\end{center}
\end{table} 

\begin{table}
	\caption{ Finding $\theta_{32,2,4}(C_{32}(R))$ and $T1_{32}(C_{32}(R))$ as $\theta_{32,2,8}(C_{32}(R))$ = $C_{32}(R)$.}
\begin{center}
\scalebox{.8}{
 \begin{tabular}{||c||c|c|c|c||}\hline \hline 
  & & & &  \\
S. No. & $C_{32}(R)$ &  $\theta_{32,2,4t}(C_{32}(R))$, $t$ = 1,2  & $C_{32}(S)\in T1_{24}(C_{32}(R))$ & T1 or T2  \\

 & & &   &  \\ \hline \hline

62.   & $C_{32}(1,4,6,15,16)$ &  = $C_{32}(4,6,7,9,16)$  & $C_{32}(1,4,6,15,16)$, $C_{32}(3,12,13,14,16)$,   &   \\ 

     & & &  $C_{32}(2,5,11,12,16)$, $C_{32}(4,7,9,10,16)$.  & T2  \\    

63.   & $C_{32}(1,4,8,10,15)$ &  = $C_{32}(4,7,8,9,10)$  & $C_{32}(1,4,8,10,15)$, $C_{32}(2,3,8,12,13)$,   &   \\ 

     & & &  $C_{32}(5,8,11,12,14)$, $C_{32}(4,6,7,8,9)$.  &  T2 \\

64.   & $C_{32}(1,4,8,12,15)$ &  = $C_{32}(4,7,8,9,12)$  & $C_{32}(1,4,8,12,15)$, $C_{32}(3,4,8,12,13)$,   &   \\
\hline \hline
\end{tabular}}
\end{center}
\end{table} 

\begin{table}
	\caption{ Finding $\theta_{32,2,4}(C_{32}(R))$ and $T1_{32}(C_{32}(R))$ as $\theta_{32,2,8}(C_{32}(R))$ = $C_{32}(R)$.}
\begin{center}
\scalebox{.8}{
 \begin{tabular}{||c||c|c|c|c||}\hline \hline 
  & & & &  \\
S. No. & $C_{32}(R)$ &  $\theta_{32,2,4t}(C_{32}(R))$, $t$ = 1,2  & $C_{32}(S)\in T1_{24}(C_{32}(R))$ & T1 or T2  \\

 & & &   &  \\ \hline \hline

93.     & $C_{32}(1,2,4,6,8,15)$ &  = $C_{32}(2,4,6,7,8,9)$  & $C_{32}(1,2,4,6,8,15)$, $C_{32}(3,6,8,12,13,14)$,   &   \\ 

     & & &  $C_{32}(2,5,8,10,11,12)$, $C_{32}(4,7,8,9,10,14)$.  &  T2 \\

94.  & $C_{32}(1,2,4,6,10,15)$ &  = $C_{32}(2,4,6,7,9,10)$  & $C_{32}(1,2,4,6,10,15)$, $C_{32}(2,3,6,12,13,14)$,   &   \\ 

     & & &  $C_{32}(2,5,10,11,12,14)$, $C_{32}(4,6,7,9,10,14)$.  & T2  \\

95.  & $C_{32}(1,2,4,6,12,15)$ &  = $C_{32}(2,4,6,7,9,12)$  & $C_{32}(1,2,4,6,12,15)$, $C_{32}(3,4,6,12,13,14)$,   &   \\
\hline \hline
\end{tabular}}
\end{center}
\end{table} 

\begin{table}
	\caption{ Finding $\theta_{32,2,4}(C_{32}(R))$ and $T1_{32}(C_{32}(R))$ as $\theta_{32,2,8}(C_{32}(R))$ = $C_{32}(R)$.}
\begin{center}
\scalebox{.78}{
 \begin{tabular}{||c||c|c|c|c||}\hline \hline 
  & & & &  \\
S. No. & $C_{32}(R)$ &  $\theta_{32,2,4t}(C_{32}(R))$, $t$ = 1,2  & $C_{32}(S)\in T1_{24}(C_{32}(R))$ & T1 or T2  \\

 & & &   &  \\ \hline \hline

124. & $C_{32}(1,2,10,12,14,15)$ & = $C_{32}(2,7,9,10,12,14)$ & $C_{32}(1,2,10,12,14,15)$, $C_{32}(2,3,4,6,10,13)$, & \\ 

     & & &  $C_{32}(4,5,6,10,11,14)$, $C_{32}(2,6,7,9,12,14)$.  &  T2 \\    

125.  & $C_{32}(1,2,10,12,15,16)$ &  = $C_{32}(2,7,9,10,12,16)$  & $C_{32}(1,2,10,12,15,16)$, $C_{32}(2,3,4,6,13,16)$,   &   \\ 

     & & &  $C_{32}(4,5,10,11,14,16)$, $C_{32}(6,7,9,12,14,16)$.  & T2  \\

126.   & $C_{32}(1,2,10,14,15,16)$ &  = $C_{32}(2,7,9,10,14,16)$  & $C_{32}(1,2,10,14,15,16)$, $C_{32}(2,3,6,10,13,16)$,   &   \\
\hline \hline
\end{tabular}}
\end{center}
\end{table} 

\begin{table}
	\caption{ Finding $\theta_{32,2,4}(C_{32}(R))$ and $T1_{32}(C_{32}(R))$ as $\theta_{32,2,8}(C_{32}(R))$ = $C_{32}(R)$.}
\begin{center}
\scalebox{.78}{
 \begin{tabular}{||c||c|c|c|c||}\hline \hline 
  & & & &  \\
S. No. & $C_{32}(R)$ &  $\theta_{32,2,4t}(C_{32}(R))$, $t$ = 1,2  & $C_{32}(S)\in T1_{24}(C_{32}(R))$ & T1 or T2  \\

 & & &   &  \\ \hline \hline

156. & $C_{32}(1,6,10,14,15,16)$ &  = $C_{32}(6,7,9,10,14,16)$  & $C_{32}(1,6,10,14,15,16)$, $C_{32}(2,3,10,13,14,16)$,   &   \\ 

     & & &  $C_{32}(2,5,6,11,14,16)$, $C_{32}(2,6,7,9,10,16)$.  &  T2 \\  

157. & $C_{32}(1,6,12,14,15,16)$ &  = $C_{32}(6,7,9,12,14,16)$  & $C_{32}(1,6,12,14,15,16)$, $C_{32}(3,4,10,13,14,16)$,   &   \\ 

     & & &  $C_{32}(2,4,5,6,11,16)$, $C_{32}(2,7,9,10,12,16)$.  & T2  \\ 

158. & $C_{32}(1,8,10,12,14,15)$ & = $C_{32}(7,8,9,10,12,14)$ & $C_{32}(1,8,10,12,14,15)$, $C_{32}(2,3,4,8,10,13)$, & \\
\hline \hline
\end{tabular}}
\end{center}
\end{table} 

\begin{table}
	\caption{ Finding $\theta_{32,2,4}(C_{32}(R))$ and $T1_{32}(C_{32}(R))$ as $\theta_{32,2,8}(C_{32}(R))$ = $C_{32}(R)$.}
\begin{center}
\scalebox{.73}{
 \begin{tabular}{||c||c|c|c|c||}\hline \hline 
  & & & &  \\
S. No. & $C_{32}(R)$ &  $\theta_{32,2,4t}(C_{32}(R))$, $t$ = 1,2  & $C_{32}(S)\in T1_{24}(C_{32}(R))$ & T1 or T2  \\

 & & &   &  \\ \hline \hline

188. & $C_{32}(1,2,6,8,14,15,16)$ &  = $C_{32}(2,6,7,8,9,14,16)$  & $C_{32}(1,2,6,8,14,15,16)$, $C_{32}(3,6,8,10,13,14,16)$, & \\ 

     & & &  $C_{32}(2,5,6,8,10,11,16)$, $C_{32}(2,7,8,9,10,14,16)$.  &  T2 \\  	

189. & $C_{32}(1,2,6,10,12,14,15)$ & = $C_{32}(2,6,7,9,10,12,14)$ & $C_{32}(1,2,6,10,12,14,15)$, $C_{32}(2,3,4,6,10,13,14)$, &  \\ 

     & & &  $C_{32}(2,4,6,5,10,11,14)$, $C_{32}(2,6,7,9,10,12,14)$.  & T1 \\ 

190. & $C_{32}(1,2,6,10,12,15,16)$ &  = $C_{32}(2,6,7,9,10,12,16)$  & $C_{32}(1,2,6,10,12,15,16)$, $C_{32}(2,3,4,6,13,14,16)$,   &   \\
\hline \hline
\end{tabular}}
\end{center}
\end{table} 

\begin{table}
	\caption{ Finding $\theta_{32,2,4}(C_{32}(R))$ and $T1_{32}(C_{32}(R))$ as $\theta_{32,2,8}(C_{32}(R))$ = $C_{32}(R)$.}
\begin{center}
\scalebox{.68}{
 \begin{tabular}{||c||c|c|c|c||}\hline \hline 
  & & & &  \\
S. No. & $C_{32}(R)$ &  $\theta_{32,2,4t}(C_{32}(R))$, $t$ = 1,2  & $C_{32}(S)\in T1_{24}(C_{32}(R))$ & T1 \\

 & & & & or \\

& & & & T2  \\

 & & &   &  \\ \hline \hline

222. & $C_{32}(1,2,4,6,8,12,14,15)$ &  = $C_{32}(2,4,6,7,8,9,12,14)$  & $C_{32}(1,2,4,6,8,12,14,15)$, $C_{32}(3,4,6,8,10,12,13,14)$,   &   \\ 

     & & &  $C_{32}(2,4,5,6,8,10,11,12)$, $C_{32}(2,4,7,8,9,10,12,14)$.  &  T2 \\  

223. & $C_{32}(1,2,4,6,8,12,15,16)$ & = $C_{32}(2,4,6,7,8,9,12,16)$ & $C_{32}(1,2,4,6,8,12,15,16)$, $C_{32}(3,4,6,8,12,13,14,16)$,   &   \\
\hline \hline
\end{tabular}}
\end{center}
\end{table} 

\vspace{.2cm}
\noindent
{\bf Case (b)}\quad  As solution to this case is similar to case (b),  to simplify our work, we present only Table 17 to Table 24. 

\begin{table}
	\caption{ Finding $\theta_{32,2,4}(C_{32}(R))$ and $T1_{32}(C_{32}(R))$ as $\theta_{32,2,8}(C_{32}(R))$ = $C_{32}(R)$.}
\begin{center}
\scalebox{.8}{
 \begin{tabular}{||c||c|c|c|c||}\hline \hline 
  & & & &  \\
S. No. & $C_{32}(R)$ &  $\theta_{32,2,4t}(C_{32}(R))$, $t$ = 1,2  & $C_{32}(S)\in T1_{24}(C_{32}(R))$ & T1 or T2  \\
 \hline \hline

1.   & $C_{32}(2,3,13)$ &  = $C_{32}(2,5,11)$  & $C_{32}(2,3,13)$, $C_{32}(6,7,9)$,   &   \\ 

     & & &  $C_{32}(1,10,15)$, $C_{32}(5,11,14)$.  &  T2 \\

2.   & $C_{32}(3,4,13)$ &  = $C_{32}(4,5,11)$  & $C_{32}(3,4,13)$, $C_{32}(7,9,12)$,   &   \\ 

     & & &  $C_{32}(1,12,15)$, $C_{32}(4,5,11)$.  &  T1 \\

3.   & $C_{32}(3,6,13)$ &  = $C_{32}(5,6,11)$  & $C_{32}(3,6,13)$, $C_{32}(7,9,14)$,    &   \\ 

     & & &  $C_{32}(1,2,15)$, $C_{32}(5,10,11)$.  &  T2 \\
\hline \hline
\end{tabular}}
\end{center}
\end{table} 

\begin{table}
	\caption{ Finding $\theta_{32,2,4}(C_{32}(R))$ and $T1_{32}(C_{32}(R))$ as $\theta_{32,2,8}(C_{32}(R))$ = $C_{32}(R)$.}
\begin{center}
\scalebox{.8}{
 \begin{tabular}{||c||c|c|c|c||}\hline \hline 
  & & & &  \\
S. No. & $C_{32}(R)$ &  $\theta_{32,2,4t}(C_{32}(R))$, $t$ = 1,2  & $C_{32}(S)\in T1_{24}(C_{32}(R))$ & T1 or T2  \\

 & & &   &  \\ \hline \hline

32.  & $C_{32}(3,10,13,14)$ &  = $C_{32}(5,10,11,14)$  & $C_{32}(3,10,13,14)$, $C_{32}(2,7,9,10)$,   &   \\ 

     & & &  $C_{32}(1,6,14,15)$, $C_{32}(2,5,6,11)$.  & T2 \\      

33.  & $C_{32}(3,10,13,16)$ &  = $C_{32}(5,10,11,16)$  & $C_{32}(3,10,13,16)$, $C_{32}(2,7,9,16)$,   &   \\ 

     & & &  $C_{32}(1,14,15,16)$, $C_{32}(5,6,11,16)$.  &  T2 \\ 

34.  & $C_{32}(3,12,13,14)$ &  = $C_{32}(5,11,12,14)$  & $C_{32}(3,12,13,14)$, $C_{32}(4,7,9,10,13)$,   &   \\
\hline \hline
\end{tabular}}
\end{center}
\end{table} 

\begin{table}
	\caption{ Finding $\theta_{32,2,4}(C_{32}(R))$ and $T1_{32}(C_{32}(R))$ as $\theta_{32,2,8}(C_{32}(R))$ = $C_{32}(R)$.}
\begin{center}
\scalebox{.8}{
 \begin{tabular}{||c||c|c|c|c||}\hline \hline 
  & & & &  \\
S. No. & $C_{32}(R)$ &  $\theta_{32,2,4t}(C_{32}(R))$, $t$ = 1,2  & $C_{32}(S)\in T1_{24}(C_{32}(R))$ & T1 or T2  \\

 & & &   &  \\ \hline \hline

62.   & $C_{32}(3,4,6,13,16)$ &  = $C_{32}(4,5,6,11,16)$  & $C_{32}(3,4,6,13,16)$, $C_{32}(7,9,12,14,16)$,   &   \\ 

     & & &  $C_{32}(1,2,12,15,16)$, $C_{32}(4,5,10,11,16)$.  & T2  \\    

63.   & $C_{32}(3,4,8,10,13)$ &  = $C_{32}(4,5,8,10,11)$  & $C_{32}(3,4,8,10,13)$, $C_{32}(2,7,8.9,12)$,   &   \\ 

     & & &  $C_{32}(1,8,12,14,15)$, $C_{32}(4,5,6,8,11)$.  &  T2 \\

64.   & $C_{32}(3,4,8,12,13)$ &  = $C_{32}(4,5,8,11,12)$  & $C_{32}(3,4,8,12,13)$, $C_{32}(4,7,8.9,12)$,   &   \\
\hline \hline
\end{tabular}}
\end{center}
\end{table} 

\begin{table}
	\caption{ Finding $\theta_{32,2,4}(C_{32}(R))$ and $T1_{32}(C_{32}(R))$ as $\theta_{32,2,8}(C_{32}(R))$ = $C_{32}(R)$.}
\begin{center}
\scalebox{.8}{
 \begin{tabular}{||c||c|c|c|c||}\hline \hline 
  & & & &  \\
S. No. & $C_{32}(R)$ &  $\theta_{32,2,4t}(C_{32}(R))$, $t$ = 1,2  & $C_{32}(S)\in T1_{24}(C_{32}(R))$ & T1 or T2  \\

 & & &   &  \\ \hline \hline

93.     & $C_{32}(2,3,4,6,8,13)$ &  = $C_{32}(2,4,5,6,8,11)$  & $C_{32}(2,3,4,6,8,13)$, $C_{32}(7,6,8,9,12,14)$,   &   \\ 

     & & &  $C_{32}(1,2,8,10,12,15)$, $C_{32}(4,5,8,10,11,14)$.  &  T2 \\

94.  & $C_{32}(2,3,4,6,10,13)$ &  = $C_{32}(2,4,5,6,10,11)$  & $C_{32}(2,3,4,6,10,13)$, $C_{32}(2,7,6,12,9,14)$,   &   \\ 

     & & &  $C_{32}(1,2,10,12,14,15)$, $C_{32}(4,6,5,10,11,14)$.  & T2  \\

95.  & $C_{32}(2,3,4,6,12,13)$ &  = $C_{32}(2,4,5,6,11,12)$  & $C_{32}(2,3,4,6,12,13)$, $C_{32}(4,6,7,9,12,13,14)$,   &   \\
\hline \hline
\end{tabular}}
\end{center}
\end{table} 

\begin{table}
	\caption{ Finding $\theta_{32,2,4}(C_{32}(R))$ and $T1_{32}(C_{32}(R))$ as $\theta_{32,2,8}(C_{32}(R))$ = $C_{32}(R)$.}
\begin{center}
\scalebox{.78}{
 \begin{tabular}{||c||c|c|c|c||}\hline \hline 
  & & & &  \\
S. No. & $C_{32}(R)$ &  $\theta_{32,2,4t}(C_{32}(R))$, $t$ = 1,2  & $C_{32}(S)\in T1_{24}(C_{32}(R))$ & T1 or T2  \\

 & & &   &  \\ \hline \hline

124. & $C_{32}(2,3,10,12,13,14)$ & = $C_{32}(2,5,10,11,12,14)$ & $C_{32}(2,3,10,12,13,14)$, $C_{32}(2,4,6,7,9,10)$, & \\ 

     & & &  $C_{32}(1,4,6,10,14,15)$, $C_{32}(2,5,6,11,12,14)$.  &  T2 \\    

125.  & $C_{32}(2,3,10,12,13,16)$ &  = $C_{32}(2,5,10,11,12,16)$  & $C_{32}(2,3,10,12,13,16)$, $C_{32}(2,4,6,7,9,16)$,   &   \\ 

     & & &  $C_{32}(1,4,10,14,15,16)$, $C_{32}(5,6,11,12,14,16)$.  & T2  \\

126.   & $C_{32}(2,3,10,13,14,16)$ &  = $C_{32}(2,5,10,11,14,16)$  & $C_{32}(2,3,10,13,14,16)$, $C_{32}(2,6,7,9,10,16)$,   &   \\
\hline \hline
\end{tabular}}
\end{center}
\end{table} 

\begin{table}
	\caption{ Finding $\theta_{32,2,4}(C_{32}(R))$ and $T1_{32}(C_{32}(R))$ as $\theta_{32,2,8}(C_{32}(R))$ = $C_{32}(R)$.}
\begin{center}
\scalebox{.78}{
 \begin{tabular}{||c||c|c|c|c||}\hline \hline 
  & & & &  \\
S. No. & $C_{32}(R)$ &  $\theta_{32,2,4t}(C_{32}(R))$, $t$ = 1,2  & $C_{32}(S)\in T1_{24}(C_{32}(R))$ & T1 or T2  \\

 & & &   &  \\ \hline \hline

156. & $C_{32}(3,6,10,14,13,16)$ &  = $C_{32}(5,6,10,11,14,16)$  & $C_{32}(3,6,10,14,13,16)$, $C_{32}(2,7,9,10,14,16)$,   &   \\ 

     & & &  $C_{32}(1,2,6,14,15,16)$, $C_{32}(2,5,6,10,11,16)$.  &  T2 \\  

157. & $C_{32}(3,6,12,14,13,16)$ &  = $C_{32}(5,6,11,12,14,16)$  & $C_{32}(3,6,12,14,13,16)$, $C_{32}(4,7,9,10,14,16)$,   &   \\ 

     & & &  $C_{32}(1,2,4,6,15,16)$, $C_{32}(2,5,10,11,12,16)$.  & T2  \\ 

158. & $C_{32}(3,8,10,12,13,14)$ & = $C_{32}(5,8,10,11,12,14)$ & $C_{32}(3,8,10,12,13,14)$, $C_{32}(2,4,7,8,9,10)$, & \\
\hline \hline
\end{tabular}}
\end{center}
\end{table} 

\begin{table}
	\caption{ Finding $\theta_{32,2,4}(C_{32}(R))$ and $T1_{32}(C_{32}(R))$ as $\theta_{32,2,8}(C_{32}(R))$ = $C_{32}(R)$.}
\begin{center}
\scalebox{.73}{
 \begin{tabular}{||c||c|c|c|c||}\hline \hline 
  & & & &  \\
S. No. & $C_{32}(R)$ &  $\theta_{32,2,4t}(C_{32}(R))$, $t$ = 1,2  & $C_{32}(S)\in T1_{24}(C_{32}(R))$ & T1 or T2  \\ 

 & & &   &  \\ \hline \hline

188. & $C_{32}(2,3,6,8,13,14,16)$ &  = $C_{32}(2,5,6,8,11,14,16)$  & $C_{32}(2,3,6,8,13,14,16)$, $C_{32}(6,7,8,9,10,14,16)$, & \\ 

     & & &  $C_{32}(1,2,6,8,10,15,16)$, $C_{32}(2,5,8,10,11,14,16)$.  &  T2 \\  	

189. & $C_{32}(2,3,6,10,12,13,14)$ & = $C_{32}(2,5,6,10,11,12,14)$ & $C_{32}(2,3,6,10,12,13,14)$, $C_{32}(2,4,6,7,9,10,14)$, &  \\ 

     & & &  $C_{32}(1,2,4,6,10,14,15)$, $C_{32}(2,5,6,10,11,12,14)$.  & T1 \\ 

190. & $C_{32}(2,3,6,10,12,13,16)$ &  = $C_{32}(2,5,6,10,11,12,16)$  & $C_{32}(2,3,6,10,12,13,16)$, $C_{32}(2,4,6,7,9,14,16)$,   &   \\
\hline \hline
\end{tabular}}
\end{center}
\end{table} 

\begin{table}
	\caption{ Finding $\theta_{32,2,4}(C_{32}(R))$ and $T1_{32}(C_{32}(R))$ as $\theta_{32,2,8}(C_{32}(R))$ = $C_{32}(R)$.}
\begin{center}
\scalebox{.68}{
 \begin{tabular}{||c||c|c|c|c||}\hline \hline 
  & & & &  \\
S. No. & $C_{32}(R)$ &  $\theta_{32,2,4t}(C_{32}(R))$, $t$ = 1,2  & $C_{32}(S)\in T1_{24}(C_{32}(R))$ & T1 \\

 & & & & or \\

& & & & T2  \\  \hline \hline

222. & $C_{32}(2,3,4,6,8,12,13,14)$ &  = $C_{32}(2,4,5,6,8,11,12,14)$  & $C_{32}(2,3,4,6,8,12,13,14)$, $C_{32}(4,6,7,8,9,10,12,14)$,   &   \\ 

     & & &  $C_{32}(1,2,4,6,8,10,12,15)$, $C_{32}(2,4,5,8,10,11,12,14)$.  &  T2 \\  

223. & $C_{32}(2,3,4,6,8,12,13,16)$ & = $C_{32}(2,4,5,6,8,11,12,16)$ & $C_{32}(2,3,4,6,8,12,13,16)$, $C_{32}(4,6,7,8,9,12,14,16)$,   &   \\ 

     & & &  $C_{32}(1,2,4,8,10,12,15,16)$, $C_{32}(4,5,8,10,11,12,14,16)$.  & T2  \\
\hline \hline
\end{tabular}}
\end{center}
\end{table} 

From the values of these tables, we could see that there are 384 pairs of isomorphic circulant graphs of order 32 and of Type-2 w.r.t. $m$ = 2.   \hfill $\Box$

All the 384 pairs of isomorphic circulant graphs of order 32 and of Type-2 w.r.t. $m$ = 2 are presented below for clarity.

\medskip
\noindent
{\bf  Pairs of isomorphic $C_{32}(R)$ and $C_{32}(S)$ of Type-2 w.r.t. $m$ = 2 $\ni$ $1,15\in R$.} 
\begin{enumerate}
\item [\rm (1)]  $C_{32} (1,2,15)$, $C_{32} (2,7,9)$;  

\item [\rm (2)] $C_{32} (1,6,15)$, $C_{32} (6,7,9)$; 
	
\item [\rm (3)] $C_{32} (1,10,15)$, $C_{32} (7,9,10)$;  
	
\item [\rm (4)] $C_{32} (1,14,15)$, $C_{32} (7,9,14)$; 
	
\item [\rm (5)] $C_{32} (1,2,4,15)$, $C_{32} (2,4,7,9)$;  
	
\item [\rm (6)] $C_{32} (1,2,6,15)$, $C_{32} (2,6,7,9)$;
	
\item [\rm (7)] $C_{32} (1,2,8,15)$, $C_{32} (2,7,8,9)$; 
	
\item [\rm (8)] $C_{32} (1,2,10,15)$, $C_{32} (2,7,9,10)$; 
	
\item [\rm (9)] $C_{32} (1,2,12,15)$, $C_{32} (2,7,9,12)$; 
	
\item [\rm (10)] $C_{32} (1,2,15,16)$, $C_{32} (2,7,9,16)$; 

\item [\rm (11)] $C_{32} (1,4,6,15)$, $C_{32} (4,6,7,9)$; 

\item [\rm (12)] $C_{32} (1,4,10,15)$, $C_{32} (4,7,9,10)$; 

\item [\rm (13)] $C_{32} (1,4,14,15)$, $C_{32} (4,7,9,14)$; 

\item [\rm (14)] $C_{32} (1,6,8,15)$, $C_{32} (6,7,8,9)$; 
	
\item [\rm (15)] $C_{32} (1,6,12,15)$, $C_{32} (6,7,9,12)$; 
	
\item [\rm (16)] $C_{32} (1,6,14,15)$, $C_{32} (6,7,9,14)$; 
	
\item [\rm (17)] $C_{32} (1,6,15,16)$, $C_{32} (6,7,9,16)$; 
	
\item [\rm (18)] $C_{32} (1,8,10,15)$, $C_{32} (7,8,9,10)$; 

\item [\rm (19)] $C_{32} (1,8,14,15)$, $C_{32} (7,8,9,14)$; 

\item [\rm (20)] $C_{32} (1,10,12,15)$, $C_{32} (7,9,10,12)$; 

\item [\rm (21)] $C_{32} (1,10,14,15)$, $C_{32} (7,9,10,14)$; 

\item [\rm (22)] $C_{32} (1,10,15,16)$, $C_{32} (7,9,10,16)$; 

\item [\rm (23)] $C_{32} (1,12,14,15)$, $C_{32} (7,9,12,14)$; 

\item [\rm (24)] $C_{32} (1,14,15,16)$, $C_{32} (7,9,14,16)$; 

\item [\rm (25)] $C_{32} (1,2,4,6,15)$, $C_{32} (2,4,6,7,9)$;

\item [\rm (26)] $C_{32} (1,2,4,8,15)$, $C_{32} (2,4,7,8,9)$;

\item [\rm (27)] $C_{32} (1,2,4,10,15)$, $C_{32} (2,4,7,9,10)$;
	
\item [\rm (28)] $C_{32} (1,2,4,12,15)$, $C_{32} (2,4,7,9,12)$; 

\item [\rm (29)] $C_{32} (1,2,4,15,16)$, $C_{32} (2,4,7,9,16)$; 

\item [\rm (30)] $C_{32} (1,2,6,8,15)$, $C_{32} (2,6,7,8,9)$;

\item [\rm (31)] $C_{32} (1,2,6,10,15)$, $C_{32} (2,6,7,9,10)$;

\item [\rm (32)] $C_{32} (1,2,6,12,15)$, $C_{32} (2,6,7,9,12)$;

\item [\rm (33)] $C_{32} (1,2,6,14,15)$, $C_{32} (2,6,7,9,14)$;

\item [\rm (34)] $C_{32} (1,2,6,15,16)$, $C_{32} (2,6,7,9,16)$;

\item [\rm (35)] $C_{32} (1,2,8,10,15)$, $C_{32} (2,7,8,9,10)$;  

\item [\rm (36)] $C_{32} (1,2,8,12,15)$, $C_{32} (2,7,8,9,12)$;  

\item [\rm (37)] $C_{32} (1,2,8,15,16)$, $C_{32} (2,7,8,9,16)$;

\item [\rm (38)] $C_{32} (1,2,10,12,15)$, $C_{32} (2,7,9,10,12)$;  

\item [\rm (39)] $C_{32} (1,2,10,14,15)$, $C_{32} (2,7,9,10,14)$;  

\item [\rm (40)] $C_{32} (1,2,10,15,16)$, $C_{32} (2,7,9,10,16)$; 

\item [\rm (41)] $C_{32} (1,2,12,15,16)$, $C_{32} (2,7,9,12,16)$;

\item [\rm (42)] $C_{32} (1,4,6,8,15)$, $C_{32} (4,6,7,8,9)$; 

\item [\rm (43)] $C_{32} (1,4,6,12,15)$, $C_{32} (4,6,7,9,12)$;

\item [\rm (44)] $C_{32} (1,4,6,14,15)$, $C_{32} (4,6,7,9,14)$;

\item [\rm (45)] $C_{32} (1,4,6,15,16)$, $C_{32} (4,6,7,9,16)$;

\item [\rm (46)] $C_{32} (1,4,8,10,15)$, $C_{32} (4,7,8,9,10)$;

\item [\rm (47)] $C_{32} (1,4,8,14,15)$, $C_{32} (4,7,8,9,14)$;

\item [\rm (48)] $C_{32} (1,4,10,12,15)$, $C_{32} (4,7,9,10,12)$;

\item [\rm (49)] $C_{32} (1,4,10,14,15)$, $C_{32} (4,7,9,10,14)$;

\item [\rm (50)] $C_{32} (1,4,10,15,16)$, $C_{32} (4,7,9,10,16)$;

\item [\rm (51)] $C_{32} (1,4,12,14,15)$, $C_{32} (4,7,9,12,14)$;

\item [\rm (52)] $C_{32} (1,4,14,15,16)$, $C_{32} (4,7,9,14,16)$;

\item [\rm (53)] $C_{32} (1,6,8,12,15)$, $C_{32} (6,7,8,9,12)$;

\item [\rm (54)] $C_{32} (1,6,8,14,15)$, $C_{32} (6,7,8,9,14)$; 

\item [\rm (55)] $C_{32} (1,6,8,15,16)$, $C_{32} (6,7,8,9,16)$; 

\item [\rm (56)] $C_{32} (1,6,10,14,15)$, $C_{32} (6,7,9,10,14)$; 

\item [\rm (57)] $C_{32} (1,6,12,14,15)$, $C_{32} (6,7,9,12,14)$; 

\item [\rm (58)] $C_{32} (1,6,12,15,16)$, $C_{32} (6,7,9,12,16)$; 

\item [\rm (59)] $C_{32} (1,6,14,15,16)$, $C_{32} (6,7,9,14,16)$; 

\item [\rm (60)] $C_{32} (1,8,10,12,15)$, $C_{32} (7,8,9,10,12)$;

\item [\rm (61)] $C_{32} (1,8,10,14,15)$, $C_{32} (7,8,9,10,14)$; 

\item [\rm (62)] $C_{32} (1,8,10,15,16)$, $C_{32} (7,8,9,10,16)$;

\item [\rm (63)] $C_{32} (1,8,12,14,15)$, $C_{32} (7,8,9,12,14)$;

\item [\rm (64)] $C_{32} (1,8,14,15,16)$, $C_{32} (7,8,9,14,16)$;

\item [\rm (65)] $C_{32} (1,10,12,14,15)$, $C_{32} (7,9,10,12,14)$; 

\item [\rm (66)] $C_{32} (1,10,12,15,16)$, $C_{32} (7,9,10,12,16)$; 

\item [\rm (67)] $C_{32} (1,10,14,15,16)$, $C_{32} (7,9,10,14,16)$; 

\item [\rm (68)] $C_{32} (1,12,14,15,16)$, $C_{32} (7,9,12,14,16)$; 

\item [\rm (69)] $C_{32} (1,2,4,6,8,15)$, $C_{32} (2,4,6,7,8,9)$; 

\item [\rm (70)] $C_{32} (1,2,4,6,10,15)$, $C_{32} (2,4,6,7,9,10)$; 

\item [\rm (71)] $C_{32} (1,2,4,6,12,15)$, $C_{32} (2,4,6,7,9,12)$; 

\item [\rm (72)] $C_{32} (1,2,4,6,14,15)$, $C_{32} (2,4,6,7,9,14)$; 

\item [\rm (73)] $C_{32} (1,2,4,6,15,16)$, $C_{32} (2,4,6,7,9,16)$; 

\item [\rm (74)] $C_{32} (1,2,4,8,10,15)$, $C_{32} (2,4,8,7,9,10)$; 

\item [\rm (75)] $C_{32} (1,2,4,8,12,15)$, $C_{32} (2,4,8,7,9,12)$; 

\item [\rm (76)] $C_{32} (1,2,4,8,15,16)$, $C_{32} (2,4,8,7,9,16)$; 

 \item [\rm (77)] $C_{32} (1,2,4,10,12,15)$, $C_{32} (2,4,7,9,10,12)$;

 \item [\rm (78)] $C_{32} (1,2,4,10,14,15)$, $C_{32} (2,4,7,9,10,14)$; 
 
 \item [\rm (79)] $C_{32} (1,2,4,10,15,16)$, $C_{32} (2,4,7,9,10,16)$; 

 \item [\rm (80)]  $C_{32} (1,2,4,12,15,16)$, $C_{32} (2,4,7,9,12,16)$; 

\item [\rm (81)]  $C_{32} (1,2,6,8,10,15)$, $C_{32} (2,6,7,8,9,10)$;

\item [\rm (82)] $C_{32} (1,2,6,8,12,15)$, $C_{32} (2,6,7,8,9,12)$; 

\item [\rm (83)] $C_{32} (1,2,6,8,14,15)$, $C_{32} (2,6,7,8,9,14)$; 

\item [\rm (84)] $C_{32} (1,2,6,8,15,16)$, $C_{32} (2,6,7,8,9,16)$; 

\item [\rm (85)] $C_{32} (1,2,6,10,12,15)$, $C_{32} (2,6,7,9,10,12)$; 

\item [\rm (86)]  $C_{32} (1,2,6,10,15,16)$, $C_{32} (2,6,7,9,10,16)$; 

\item [\rm (87)] $C_{32} (1,2,6,12,14,15)$, $C_{32} (2,6,7,9,12,14)$; 

\item [\rm (88)] $C_{32} (1,2,6,12,15,16)$, $C_{32} (2,6,7,9,12,16)$;

\item [\rm (89)] $C_{32} (1,2,6,14,15,16)$, $C_{32} (2,6,7,9,14,16)$; 

\item [\rm (90)] $C_{32} (1,2,8,10,12,15)$, $C_{32} (2,7,8,9,10,12)$; 

\item [\rm (91)] $C_{32} (1,2,8,10,14,15)$, $C_{32} (2,7,8,9,10,14)$; 

\item [\rm (92)] $C_{32} (1,2,8,10,15,16)$, $C_{32} (2,7,8,9,10,16)$; 

\item [\rm (93)] $C_{32} (1,2,8,12,15,16)$, $C_{32} (2,7,8,9,12,16)$; 

\item [\rm (94)]  $C_{32} (1,2,10,12,14,15)$, $C_{32} (2,7,9,10,12,14)$; 

\item [\rm (95)] $C_{32} (1,2,10,12,15,16)$, $C_{32} (2,7,9,10,12,16)$; 

\item [\rm (96)] $C_{32} (1,2,10,14,15,16)$, $C_{32} (2,7,9,10,14,16)$; 
 
\item [\rm (97)] $C_{32} (1,4,6,8,12,15)$, $C_{32} (4,6,7,8,9,12)$; 

\item [\rm (98)] $C_{32} (1,4,6,8,14,15)$, $C_{32} (4,6,7,8,9,14)$; 

\item [\rm (99)] $C_{32} (1,4,6,8,15,16)$, $C_{32} (4,6,7,8,9,16)$; 

\item [\rm (100)] $C_{32} (1,4,6,10,14,15)$, $C_{32} (4,6,7,9,10,14)$; 

 \item [\rm (101)] $C_{32} (1,4,6,12,14,15)$, $C_{32} (4,6,7,9,12,14)$; 

\item [\rm (102)] $C_{32} (1,4,6,12,15,16)$, $C_{32} (4,6,7,9,12,16)$; 

\item [\rm (103)] $C_{32} (1,4,6,14,15,16)$, $C_{32} (4,6,7,9,14,16)$; 

\item [\rm (104)] $C_{32} (1,4,8,10,12,15)$, $C_{32} (4,7,8,9,10,12)$; 

\item [\rm (105)] $C_{32} (1,4,8,10,14,15)$, $C_{32} (4,7,8,9,10,14)$; 

\item [\rm (106)] $C_{32} (1,4,8,10,15,16)$, $C_{32} (4,7,8,9,10,16)$; 

\item [\rm (107)] $C_{32} (1,4,8,12,14,15)$, $C_{32} (4,7,8,9,12,14)$; 

\item [\rm (108)] $C_{32} (1,4,8,14,15,16)$, $C_{32} (4,7,8,9,14,16)$; 

\item [\rm (109)] $C_{32} (1,4,10,12,14,15)$, $C_{32} (4,7,9,10,12,14)$; 

\item [\rm (110)] $C_{32} (1,4,10,12,15,16)$, $C_{32} (4,7,9,10,12,16)$; 

\item [\rm (111)] $C_{32} (1,4,10,14,15,16)$, $C_{32} (4,7,9,10,14,16)$; 

\item [\rm (112)] $C_{32} (1,4,12,14,15,16)$, $C_{32} (4,7,9,12,14,16)$; 

\item [\rm (113)] $C_{32} (1,6,8,10,14,15)$, $C_{32} (6,7,8,9,10,14)$; 

 \item [\rm (114)] $C_{32} (1,6,8,12,14,15)$, $C_{32} (6,7,8,9,12,14)$; 

 \item [\rm (115)] $C_{32} (1,6,8,12,15,16)$, $C_{32} (6,7,8,9,12,16)$; 

\item [\rm (116)] $C_{32} (1,6,8,14,15,16)$, $C_{32} (6,7,8,9,14,16)$; 

 \item [\rm (117)] $C_{32} (1,6,10,12,14,15)$, $C_{32} (6,7,9,10,12,14)$; 

 \item [\rm (118)] $C_{32} (1,6,10,14,15,16)$, $C_{32} (6,7,9,10,14,16)$; 

 \item [\rm (119)] $C_{32} (1,6,12,14,15,16)$, $C_{32} (6,7,9,12,14,16)$; 

 \item [\rm (120)] $C_{32} (1,8,10,12,14,15)$, $C_{32} (7,8,9,10,12,14)$; 

 \item [\rm (121)] $C_{32} (1,8,10,12,15,16)$, $C_{32} (7,8,9,10,12,16)$; 

 \item [\rm (122)] $C_{32} (1,8,10,14,15,16)$, $C_{32} (7,8,9,10,14,16)$; 

 \item [\rm (123)] $C_{32} (1,8,12,14,15,16)$, $C_{32} (7,8,9,12,14,16)$; 

 \item [\rm (124)] $C_{32} (1,10,12,14,15,16)$, $C_{32} (7,9,10,12,14,16)$; 

 \item [\rm (125)] $C_{32} (1,2,4,6,8,10,15)$, $C_{32} (2,4,6,7,8,9,10)$; 

 \item [\rm (126)] $C_{32} (1,2,4,6,8,12,15)$, $C_{32} (2,4,6,7,8,9,12)$; 

 \item [\rm (127)] $C_{32} (1,2,4,6,8,14,15)$, $C_{32} (2,4,6,7,8,9,14)$; 

 \item [\rm (128)] $C_{32} (1,2,4,6,8,15,16)$, $C_{32} (2,4,6,7,8,9,16)$; 

  \item [\rm (129)] $C_{32} (1,2,4,6,10,12,15)$, $C_{32} (2,4,6,7,10,9,12)$; 

 \item [\rm (130)] $C_{32} (1,2,4,6,10,15,16)$, $C_{32} (2,4,6,7,9,10,16)$; 

 \item [\rm (131)] $C_{32} (1,2,4,6,12,14,15)$, $C_{32} (2,4,6,7,9,12,14)$; 

 \item [\rm (132)] $C_{32} (1,2,4,6,12,15,16)$, $C_{32} (2,4,6,7,9,12,16)$; 

 \item [\rm (133)] $C_{32} (1,2,4,6,14,15,16)$, $C_{32} (2,4,6,7,9,14,16)$; 

 \item [\rm (134)] $C_{32} (1,2,4,8,10,12,15)$, $C_{32} (2,4,7,8,9,10,12)$; 

 \item [\rm (135)] $C_{32} (1,2,4,8,10,14,15)$, $C_{32} (2,4,7,8,9,10,14)$; 

 \item [\rm (136)] $C_{32} (1,2,4,8,10,15,16)$, $C_{32} (2,4,7,8,9,10,16)$; 

 \item [\rm (137)] $C_{32} (1,2,4,8,12,15,16)$, $C_{32} (2,4,7,8,9,12,16)$; 

 \item [\rm (138)] $C_{32} (1,2,4,10,12,14,15)$, $C_{32} (2,4,7,9,10,12,14)$; 

 \item [\rm (139)] $C_{32} (1,2,4,10,12,15,16)$, $C_{32} (2,4,7,9,10,12,16)$; 

 \item [\rm (140)] $C_{32} (1,2,4,10,14,15,16)$, $C_{32} (2,4,7,9,10,14,16)$; 

 \item [\rm (141)] $C_{32} (1,2,6,8,10,12,15)$, $C_{32} (2,6,7,8,9,10,12)$; 

 \item [\rm (142)] $C_{32} (1,2,6,8,10,15,16)$, $C_{32} (2,6,7,8,9,10,16)$; 

 \item [\rm (143)] $C_{32} (1,2,6,8,12,14,15)$, $C_{32} (2,6,7,8,9,12,14)$; 

\item [\rm (144)] $C_{32} (1,2,6,8,12,15,16)$, $C_{32} (2,6,7,8,9,12,16)$; 

\item [\rm (145)] $C_{32} (1,2,6,8,14,15,16)$, $C_{32} (2,6,7,8,9,14,16)$; 

\item [\rm (146)] $C_{32} (1,2,6,10,12,15,16)$, $C_{32} (2,6,7,9,10,12,16)$; 

\item [\rm (147)] $C_{32} (1,2,6,12,14,15,16)$, $C_{32} (2,6,7,9,12,14,16)$; 

\item [\rm (148)] $C_{32} (1,2,8,10,12,14,15)$, $C_{32} (2,7,8,9,10,12,14)$; 

\item [\rm (149)] $C_{32} (1,2,8,10,12,15,16)$, $C_{32} (2,7,8,9,10,12,16)$; 

\item [\rm (150)] $C_{32} (1,2,8,10,14,15,16)$, $C_{32} (2,7,8,9,10,14,16)$; 

\item [\rm (151)] $C_{32} (1,2,10,12,14,15,16)$, $C_{32} (2,7,9,10,12,14,16)$; 

\item [\rm (152)] $C_{32} (1,4,6,8,10,14,15)$, $C_{32} (4,6,7,8,9,10,14)$; 

\item [\rm (153)] $C_{32} (1,4,6,8,12,14,15)$, $C_{32} (4,6,7,8,9,12,14)$; 

\item [\rm (154)] $C_{32} (1,4,6,8,12,15,16)$, $C_{32} (4,6,7,8,9,12,16)$; 

\item [\rm (155)] $C_{32} (1,4,6,8,14,15,16)$, $C_{32} (4,6,7,8,9,14,16)$; 

\item [\rm (156)] $C_{32} (1,4,6,10,12,14,15)$, $C_{32} (4,6,7,9,10,12,14)$; 

\item [\rm (157)] $C_{32} (1,4,6,10,14,15,16)$, $C_{32} (4,6,7,9,10,14,16)$; 

\item [\rm (158)] $C_{32} (1,4,6,12,14,15,16)$, $C_{32} (4,6,7,9,12,14,16)$; 

\item [\rm (159)] $C_{32} (1,4,8,10,12,14,15)$, $C_{32} (4,7,8,9,10,12,14)$; 

\item [\rm (160)] $C_{32} (1,4,8,10,12,15,16)$, $C_{32} (4,7,8,9,10,12,16)$; 

\item [\rm (161)] $C_{32} (1,4,8,10,14,15,16)$, $C_{32} (4,7,8,9,10,14,16)$; 

\item [\rm (162)] $C_{32} (1,4,8,12,14,15,16)$, $C_{32} (4,7,8,9,12,14,16)$; 

\item [\rm (163)] $C_{32} (1,4,10,12,14,15,16)$, $C_{32} (4,7,9,10,12,14,16)$; 

\item [\rm (164)] $C_{32} (1,6,8,10,12,14,15)$, $C_{32} (6,7,8,9,10,12,14)$; 

\item [\rm (165)] $C_{32} (1,6,8,10,14,15,16)$, $C_{32} (6,7,8,9,10,14,16)$; 

\item [\rm (166)] $C_{32} (1,6,8,12,14,15,16)$, $C_{32} (6,7,8,9,12,14,16)$; 

\item [\rm (167)] $C_{32} (1,6,10,12,14,15,16)$, $C_{32} (6,7,9,10,12,14,16)$;

\item [\rm (168)] $C_{32} (1,8,10,12,14,15,16)$, $C_{32} (7,8,9,10,12,14,16)$; 
 
\item [\rm (169)] $C_{32} (1,2,4,6,8,10,12,15)$, $C_{32} (2,4,6,7,8,9,10,12)$; 

\item [\rm (170)] $C_{32} (1,2,4,6,8,10,15,16)$, $C_{32} (2,4,6,7,8,9,10,16)$;

\item [\rm (171)] $C_{32} (1,2,4,6,8,12,14,15)$, $C_{32} (2,4,6,7,8,9,12,14)$; 

\item [\rm (172)] $C_{32} (1,2,4,6,8,12,15,16)$, $C_{32} (2,4,6,7,8,9,12,16)$; 

\item [\rm (173)] $C_{32} (1,2,4,6,8,14,15,16)$, $C_{32} (2,4,6,7,8,9,14,16)$; 

 \item [\rm (174)] $C_{32} (1,2,4,6,10,12,15,16)$, $C_{32} (2,4,6,7,9,10,12,16)$; 

 \item [\rm (175)] $C_{32} (1,2,4,6,12,14,15,16)$, $C_{32} (2,4,6,7,9,12,14,16)$; 

 \item [\rm (176)] $C_{32} (1,2,4,8,10,12,14,15)$, $C_{32} (2,4,7,8,9,10,12,14)$; 

\item [\rm (177)] $C_{32} (1,2,4,8,10,12,15,16)$, $C_{32} (2,4,7,8,9,10,12,16)$; 

\item [\rm (178)] $C_{32} (1,2,4,8,10,14,15,16)$, $C_{32} (2,4,7,8,9,10,14,16)$; 

\item [\rm (179)] $C_{32} (1,2,4,10,12,14,15,16)$, $C_{32} (2,4,7,9,10,12,14,16)$; 

\item [\rm (180)] $C_{32} (1,2,6,8,10,12,15,16)$, $C_{32} (2,6,7,8,9,10,12,16)$; 

\item [\rm (181)] $C_{32} (1,2,6,8,12,14,15,16)$, $C_{32} (2,6,7,8,9,12,14,16)$; 

\item [\rm (182)] $C_{32} (1,2,8,10,12,14,15,16)$, $C_{32} (2,7,8,9,10,12,14,16)$; 

\item [\rm (183)] $C_{32} (1,4,6,8,10,12,14,15)$, $C_{32} (4,6,7,8,9,10,12,14)$; 

\item [\rm (184)] $C_{32} (1,4,6,8,10,14,15,16)$, $C_{32} (4,6,7,8,9,10,14,16)$; 

\item [\rm (185)] $C_{32} (1,4,6,8,12,14,15,16)$, $C_{32} (4,6,7,8,9,12,14,16)$; 

\item [\rm (186)] $C_{32} (1,4,6,10,12,14,15,16)$, $C_{32} (4,6,7,9,10,12,14,16)$; 

\item [\rm (187)] $C_{32} (1,4,8,10,12,14,15,16)$, $C_{32} (4,7,8,9,10,12,14,16)$; 

\item [\rm (188)] $C_{32} (1,6,8,10,12,14,15,16)$, $C_{32} (6,7,8,9,10,12,14,16)$; 

\item [\rm (189)] $C_{32} (1,2,4,6,8,10,12,15,16)$, $C_{32} (2,4,6,7,8,9,10,12,16)$; 

\item [\rm (190)] $C_{32} (1,2,4,6,8,12,14,15,16)$, $C_{32} (2,4,6,7,8,9,12,14,16)$; 

\item [\rm (191)] $C_{32} (1,2,4,8,10,12,14,15,16)$, $C_{32} (2,4,7,8,9,10,12,14,16)$; 

\item [\rm (192)] $C_{32} (1,4,6,8,10,12,14,15,16)$, $C_{32} (4,6,7,8,9,10,12,14,16)$. 
\end{enumerate} 

\noindent
\vspace{.2cm}
{\bf  Pairs of isomorphic $C_{32}(R)$ and $C_{32}(S)$ of Type-2 w.r.t. $m$ = 2 $\ni$ $3,13\in R$.} 
\begin{enumerate}
\item [\rm (1)]  $C_{32} (2,3,13)$, $C_{32} (2,5,11)$;  

\item [\rm (2)] $C_{32} (3,6,13)$, $C_{32} (5,6,11)$; 
	
\item [\rm (3)] $C_{32} (3,10,13)$, $C_{32} (5,10,11)$;  
	
\item [\rm (4)] $C_{32} (3,13,14)$, $C_{32} (5,11,14)$; 
	
\item [\rm (5)] $C_{32} (2,3,4,13)$, $C_{32} (2,4,5,11)$;  
	
\item [\rm (6)] $C_{32} (2,3,6,13)$, $C_{32} (2,5,6,11)$;
	
\item [\rm (7)] $C_{32} (2,3,8,13)$, $C_{32} (2,5,8,11)$; 
	
\item [\rm (8)] $C_{32} (2,3,10,13)$, $C_{32} (2,5,10,11)$; 
	
\item [\rm (9)] $C_{32} (2,3,12,13)$, $C_{32} (2,5,11,12)$; 
	
\item [\rm (10)] $C_{32} (2,3,13,16)$, $C_{32} (2,5,11,16)$; 

\item [\rm (11)] $C_{32} (3,4,6,13)$, $C_{32} (4,5,6,11)$; 

\item [\rm (12)] $C_{32} (3,4,10,13)$, $C_{32} (4,5,10,11)$; 

\item [\rm (13)] $C_{32} (3,4,13,14)$, $C_{32} (4,5,11,14)$; 

\item [\rm (14)] $C_{32} (3,6,8,13)$, $C_{32} (5,6,8,11)$; 
	
\item [\rm (15)] $C_{32} (3,6,12,13)$, $C_{32} (5,6,11,12)$; 
	
\item [\rm (16)] $C_{32} (3,6,13,14)$, $C_{32} (5,6,11,14)$; 
	
\item [\rm (17)] $C_{32} (3,6,13,16)$, $C_{32} (5,6,11,16)$; 
	
\item [\rm (18)] $C_{32} (3,8,10,13)$, $C_{32} (5,8,10,11)$; 

\item [\rm (19)] $C_{32} (3,8,13,14)$, $C_{32} (5,8,11,14)$; 

\item [\rm (20)] $C_{32} (3,10,12,13)$, $C_{32} (5,10,11,12)$; 

\item [\rm (21)] $C_{32} (3,10,13,14)$, $C_{32} (5,10,11,14)$; 

\item [\rm (22)] $C_{32} (3,10,13,16)$, $C_{32} (5,10,11,16)$; 

\item [\rm (23)] $C_{32} (3,12,13,14)$, $C_{32} (5,11,12,14)$; 

\item [\rm (24)] $C_{32} (3,13,14,16)$, $C_{32} (5,11,14,16)$; 

\item [\rm (25)] $C_{32} (2,3,4,6,13)$, $C_{32} (2,4,5,6,11)$;

\item [\rm (26)] $C_{32} (2,3,4,8,13)$, $C_{32} (2,4,5,8,11)$;

\item [\rm (27)] $C_{32} (2,3,4,10,13)$, $C_{32} (2,4,5,10,11)$;
	
\item [\rm (28)] $C_{32} (2,3,4,12,13)$, $C_{32} (2,4,5,11,12)$; 

\item [\rm (29)] $C_{32} (2,3,4,13,16)$, $C_{32} (2,4,5,11,16)$; 

\item [\rm (30)] $C_{32} (2,3,6,8,13)$, $C_{32} (2,5,6,8,11)$;

\item [\rm (31)] $C_{32} (2,3,6,10,13)$, $C_{32} (2,5,6,10,11)$;

\item [\rm (32)] $C_{32} (2,3,6,12,13)$, $C_{32} (2,5,6,11,12)$;

\item [\rm (33)] $C_{32} (2,3,6,13,14)$, $C_{32} (2,5,6,11,14)$;

\item [\rm (34)] $C_{32} (2,3,6,13,16)$, $C_{32} (2,5,6,11,16)$;

\item [\rm (35)] $C_{32} (2,3,8,10,13)$, $C_{32} (2,5,8,10,11)$;  

\item [\rm (36)] $C_{32} (2,3,8,12,13)$, $C_{32} (2,5,8,11,12)$;  

\item [\rm (37)] $C_{32} (2,3,8,13,16)$, $C_{32} (2,5,8,11,16)$;

\item [\rm (38)] $C_{32} (2,3,10,12,13)$, $C_{32} (2,5,10,11,12)$;  

\item [\rm (39)] $C_{32} (2,3,10,13,14)$, $C_{32} (2,5,10,11,14)$;  

\item [\rm (40)] $C_{32} (2,3,10,13,16)$, $C_{32} (2,5,10,11,16)$; 

\item [\rm (41)] $C_{32} (2,3,12,13,16)$, $C_{32} (2,5,11,12,16)$;

\item [\rm (42)] $C_{32} (3,4,6,8,13)$, $C_{32} (4,5,6,8,11)$; 

\item [\rm (43)] $C_{32} (3,4,6,12,13)$, $C_{32} (4,5,6,11,12)$;

\item [\rm (44)] $C_{32} (3,4,6,13,14)$, $C_{32} (4,5,6,11,14)$;

\item [\rm (45)] $C_{32} (3,4,6,13,16)$, $C_{32} (4,5,6,11,16)$;

\item [\rm (46)] $C_{32} (3,4,8,10,13)$, $C_{32} (4,5,8,10,11)$;

\item [\rm (47)] $C_{32} (3,4,8,13,14)$, $C_{32} (4,5,8,11,14)$;

\item [\rm (48)] $C_{32} (3,4,10,12,13)$, $C_{32} (4,5,10,11,12)$;

\item [\rm (49)] $C_{32} (3,4,10,13,14)$, $C_{32} (4,5,10,11,14)$;

\item [\rm (50)] $C_{32} (3,4,10,13,16)$, $C_{32} (4,5,10,11,16)$;

\item [\rm (51)] $C_{32} (3,4,12,13,14)$, $C_{32} (4,5,11,12,14)$;

\item [\rm (52)] $C_{32} (3,4,13,14,16)$, $C_{32} (4,5,11,14,16)$;

\item [\rm (53)] $C_{32} (3,6,8,12,13)$, $C_{32} (5,6,8,11,12)$;

\item [\rm (54)] $C_{32} (3,6,8,13,14)$, $C_{32} (5,6,8,11,14)$; 

\item [\rm (55)] $C_{32} (3,6,8,13,16)$, $C_{32} (5,6,8,11,16)$; 

\item [\rm (56)] $C_{32} (3,6,10,13,14)$, $C_{32} (5,6,10,11,14)$; 

\item [\rm (57)] $C_{32} (3,6,12,13,14)$, $C_{32} (5,6,11,12,14)$; 

\item [\rm (58)] $C_{32} (3,6,12,13,16)$, $C_{32} (5,6,11,12,16)$; 

\item [\rm (59)] $C_{32} (3,6,13,14,16)$, $C_{32} (5,6,11,14,16)$; 

\item [\rm (60)] $C_{32} (3,8,10,12,13)$, $C_{32} (5,8,10,11,12)$;

\item [\rm (61)] $C_{32} (3,8,10,13,14)$, $C_{32} (5,8,10,11,14)$; 

\item [\rm (62)] $C_{32} (3,8,10,13,16)$, $C_{32} (5,8,10,11,16)$;

\item [\rm (63)] $C_{32} (3,8,12,13,14)$, $C_{32} (5,8,11,12,14)$;

\item [\rm (64)] $C_{32} (3,8,13,14,16)$, $C_{32} (5,8,11,14,16)$;

\item [\rm (65)] $C_{32} (3,10,12,13,14)$, $C_{32} (5,10,11,12,14)$; 

\item [\rm (66)] $C_{32} (3,10,12,13,16)$, $C_{32} (5,10,11,12,16)$; 

\item [\rm (67)] $C_{32} (3,10,13,14,16)$, $C_{32} (5,10,11,14,16)$; 

\item [\rm (68)] $C_{32} (3,12,13,14,16)$, $C_{32} (5,11,12,14,16)$; 

\item [\rm (69)] $C_{32} (2,3,4,6,8,13)$, $C_{32} (2,4,5,6,8,11)$; 

\item [\rm (70)] $C_{32} (2,3,4,6,10,13)$, $C_{32} (2,4,5,6,10,11)$; 

\item [\rm (71)] $C_{32} (2,3,4,6,12,13)$, $C_{32} (2,4,5,6,11,12)$; 

\item [\rm (72)] $C_{32} (2,3,4,6,13,14)$, $C_{32} (2,4,5,6,11,14)$; 

\item [\rm (73)] $C_{32} (2,3,4,6,13,16)$, $C_{32} (2,4,5,6,11,16)$; 

\item [\rm (74)] $C_{32} (2,3,4,8,10,13)$, $C_{32} (2,4,8,5,10,11)$; 

\item [\rm (75)] $C_{32} (2,3,4,8,12,13)$, $C_{32} (2,4,8,5,11,12)$; 

\item [\rm (76)] $C_{32} (2,3,4,8,13,16)$, $C_{32} (2,4,8,5,11,16)$; 

 \item [\rm (77)] $C_{32} (2,3,4,10,12,13)$, $C_{32} (2,4,5,10,11,12)$;

 \item [\rm (78)] $C_{32} (2,3,4,10,13,14)$, $C_{32} (2,4,5,10,11,14)$; 
 
 \item [\rm (79)] $C_{32} (2,3,4,10,13,16)$, $C_{32} (2,4,5,10,11,16)$; 

\item [\rm (80)] $C_{32} (2,3,4,12,13,16)$, $C_{32} (2,4,5,11,12,16)$; 

\item [\rm (81)] $C_{32} (2,3,6,8,10,13)$, $C_{32} (2,5,6,8,10,11)$;

\item [\rm (82)] $C_{32} (2,3,6,8,12,13)$, $C_{32} (2,5,6,8,11,12)$; 

\item [\rm (83)] $C_{32} (2,3,6,8,13,14)$, $C_{32} (2,5,6,8,11,14)$; 

\item [\rm (84)] $C_{32} (2,3,6,8,13,16)$, $C_{32} (2,5,6,8,11,16)$; 

\item [\rm (85)] $C_{32} (2,3,6,10,12,13)$, $C_{32} (2,5,6,10,11,12)$; 

\item [\rm (86)] $C_{32} (2,3,6,10,13,16)$, $C_{32} (2,5,6,10,11,16)$; 

\item [\rm (87)] $C_{32} (2,3,6,12,13,14)$, $C_{32} (2,5,6,11,12,14)$; 

\item [\rm (88)] $C_{32} (2,3,6,12,13,16)$, $C_{32} (2,5,6,11,12,16)$;

\item [\rm (89)] $C_{32} (2,3,6,13,14,16)$, $C_{32} (2,5,6,11,14,16)$; 

\item [\rm (90)] $C_{32} (2,3,8,10,12,13)$, $C_{32} (2,5,8,10,11,12)$; 

\item [\rm (91)] $C_{32} (2,3,8,10,13,14)$, $C_{32} (2,5,8,10,11,14)$; 

\item [\rm (92)] $C_{32} (2,3,8,10,13,16)$, $C_{32} (2,5,8,10,11,16)$; 

\item [\rm (93)] $C_{32} (2,3,8,12,13,16)$, $C_{32} (2,5,8,11,12,16)$; 

\item [\rm (94)] $C_{32} (2,3,10,12,13,14)$, $C_{32} (2,5,10,11,12,14)$; 

\item [\rm (95)] $C_{32} (2,3,10,12,13,16)$, $C_{32} (2,5,10,11,12,16)$; 

\item [\rm (96)] $C_{32} (2,3,10,13,14,16)$, $C_{32} (2,5,10,11,14,16)$; 
 
\item [\rm (97)] $C_{32} (3,4,6,8,12,13)$, $C_{32} (4,5,6,8,11,12)$; 

\item [\rm (98)] $C_{32} (3,4,6,8,13,14)$, $C_{32} (4,5,6,8,11,14)$; 

\item [\rm (99)] $C_{32} (3,4,6,8,13,16)$, $C_{32} (4,5,6,8,11,16)$; 

\item [\rm (100)] $C_{32} (3,4,6,10,13,14)$, $C_{32} (4,5,6,10,11,14)$; 

\item [\rm (101)] $C_{32} (3,4,6,12,13,14)$, $C_{32} (4,5,6,11,12,14)$; 

\item [\rm (102)] $C_{32} (3,4,6,12,13,16)$, $C_{32} (4,5,6,11,12,16)$; 

\item [\rm (103)] $C_{32} (3,4,6,13,14,16)$, $C_{32} (4,5,6,11,14,16)$; 

\item [\rm (104)] $C_{32} (3,4,8,10,12,13)$, $C_{32} (4,5,8,10,11,12)$; 

\item [\rm (105)] $C_{32} (3,4,8,10,13,14)$, $C_{32} (4,5,8,10,11,14)$; 

\item [\rm (106)] $C_{32} (3,4,8,10,13,16)$, $C_{32} (4,5,8,10,11,16)$; 

\item [\rm (107)] $C_{32} (3,4,8,12,13,14)$, $C_{32} (4,5,8,11,12,14)$; 

\item [\rm (108)] $C_{32} (3,4,8,13,14,16)$, $C_{32} (4,5,8,11,14,16)$; 

\item [\rm (109)] $C_{32} (3,4,10,12,13,14)$, $C_{32} (4,5,10,11,12,14)$; 

\item [\rm (110)] $C_{32} (3,4,10,12,13,16)$, $C_{32} (4,5,10,11,12,16)$; 

\item [\rm (111)] $C_{32} (3,4,10,13,14,16)$, $C_{32} (4,5,10,11,14,16)$; 

\item [\rm (112)] $C_{32} (3,4,12,13,14,16)$, $C_{32} (4,5,11,12,14,16)$; 

\item [\rm (113)] $C_{32} (3,6,8,10,13,14)$, $C_{32} (5,6,8,10,11,14)$; 

\item [\rm (114)] $C_{32} (3,6,8,12,13,14)$, $C_{32} (5,6,8,11,12,14)$; 

 \item [\rm (115)] $C_{32} (3,6,8,12,13,16)$, $C_{32} (5,6,8,11,12,16)$; 

\item [\rm (116)] $C_{32} (3,6,8,13,14,16)$, $C_{32} (5,6,8,11,14,16)$; 

 \item [\rm (117)] $C_{32} (3,6,10,12,13,14)$, $C_{32} (5,6,10,11,12,14)$; 

 \item [\rm (118)] $C_{32} (3,6,10,13,14,16)$, $C_{32} (5,6,10,11,14,16)$; 

 \item [\rm (119)] $C_{32} (3,6,12,13,14,16)$, $C_{32} (5,6,11,12,14,16)$; 

 \item [\rm (120)] $C_{32} (3,8,10,12,13,14)$, $C_{32} (5,8,10,11,12,14)$; 

 \item [\rm (121)] $C_{32} (3,8,10,12,13,16)$, $C_{32} (5,8,10,11,12,16)$; 

 \item [\rm (122)] $C_{32} (3,8,10,13,14,16)$, $C_{32} (5,8,10,11,14,16)$; 

 \item [\rm (123)] $C_{32} (3,8,12,13,14,16)$, $C_{32} (5,8,11,12,14,16)$; 

 \item [\rm (124)] $C_{32} (3,10,12,13,14,16)$, $C_{32} (5,10,11,12,14,16)$; 

 \item [\rm (125)] $C_{32} (2,3,4,6,8,10,13)$, $C_{32} (2,4,5,6,8,10,11)$; 

 \item [\rm (126)] $C_{32} (2,3,4,6,8,12,13)$, $C_{32} (2,4,5,6,8,11,12)$; 

 \item [\rm (127)] $C_{32} (2,3,4,6,8,13,14)$, $C_{32} (2,4,5,6,8,11,14)$; 

 \item [\rm (128)] $C_{32} (2,3,4,6,8,13,16)$, $C_{32} (2,4,5,6,8,11,16)$; 

  \item [\rm (129)] $C_{32} (2,3,4,6,10,12,13)$, $C_{32} (2,4,5,6,10,11,12)$; 

 \item [\rm (130)] $C_{32} (2,3,4,6,10,13,16)$, $C_{32} (2,4,5,6,10,11,16)$; 

 \item [\rm (131)] $C_{32} (2,3,4,6,12,13,14)$, $C_{32} (2,4,5,6,11,12,14)$; 

 \item [\rm (132)] $C_{32} (2,3,4,6,12,13,16)$, $C_{32} (2,4,5,6,11,12,16)$; 

 \item [\rm (133)] $C_{32} (2,3,4,6,13,14,16)$, $C_{32} (2,4,5,6,11,14,16)$; 

 \item [\rm (134)] $C_{32} (2,3,4,8,10,12,13)$, $C_{32} (2,4,5,8,10,11,12)$; 

 \item [\rm (135)] $C_{32} (2,3,4,8,10,13,14)$, $C_{32} (2,4,5,8,10,11,14)$; 

 \item [\rm (136)] $C_{32} (2,3,4,8,10,13,16)$, $C_{32} (2,4,5,8,10,11,16)$; 

 \item [\rm (137)] $C_{32} (2,3,4,8,12,13,16)$, $C_{32} (2,4,5,8,11,12,16)$; 

 \item [\rm (138)] $C_{32} (2,3,4,10,12,13,14)$, $C_{32} (2,4,5,10,11,12,14)$; 

 \item [\rm (139)] $C_{32} (2,3,4,10,12,13,16)$, $C_{32} (2,4,5,10,11,12,16)$; 

 \item [\rm (140)] $C_{32} (2,3,4,10,13,14,16)$, $C_{32} (2,4,5,10,11,14,16)$; 

 \item [\rm (141)] $C_{32} (2,3,6,8,10,12,13)$, $C_{32} (2,5,6,8,10,11,12)$; 

 \item [\rm (142)] $C_{32} (2,3,6,8,10,13,16)$, $C_{32} (2,5,6,8,10,11,16)$; 

 \item [\rm (143)] $C_{32} (2,3,6,8,12,13,14)$, $C_{32} (2,5,6,8,11,12,14)$; 

\item [\rm (144)] $C_{32} (2,3,6,8,12,13,16)$, $C_{32} (2,5,6,8,11,12,16)$; 

\item [\rm (145)] $C_{32} (2,3,6,8,13,14,16)$, $C_{32} (2,5,6,8,11,14,16)$; 

\item [\rm (146)] $C_{32} (2,3,6,10,12,13,16)$, $C_{32} (2,5,6,10,11,12,16)$; 

\item [\rm (147)] $C_{32} (2,3,6,12,13,14,16)$, $C_{32} (2,5,6,11,12,14,16)$; 

\item [\rm (148)] $C_{32} (2,3,8,10,12,13,14)$, $C_{32} (2,5,8,10,11,12,14)$; 

\item [\rm (149)] $C_{32} (2,3,8,10,12,13,16)$, $C_{32} (2,5,8,10,11,12,16)$; 

\item [\rm (150)] $C_{32} (2,3,8,10,13,14,16)$, $C_{32} (2,5,8,10,11,14,16)$; 

\item [\rm (151)] $C_{32} (2,3,10,12,13,14,16)$, $C_{32} (2,5,10,11,12,14,16)$; 

\item [\rm (152)] $C_{32} (3,4,6,8,10,13,14)$, $C_{32} (4,5,6,8,10,11,14)$; 

\item [\rm (153)] $C_{32} (3,4,6,8,12,13,14)$, $C_{32} (4,5,6,8,11,12,14)$; 

\item [\rm (154)] $C_{32} (3,4,6,8,12,13,16)$, $C_{32} (4,5,6,8,11,12,16)$; 

\item [\rm (155)] $C_{32} (3,4,6,8,13,14,16)$, $C_{32} (4,5,6,8,11,14,16)$; 

\item [\rm (156)] $C_{32} (3,4,6,10,12,13,14)$, $C_{32} (4,5,6,10,11,12,14)$; 

\item [\rm (157)] $C_{32} (3,4,6,10,13,14,16)$, $C_{32} (4,5,6,10,11,14,16)$; 

\item [\rm (158)] $C_{32} (3,4,6,12,13,14,16)$, $C_{32} (4,5,6,11,12,14,16)$; 

\item [\rm (159)] $C_{32} (3,4,8,10,12,13,14)$, $C_{32} (4,5,8,10,11,12,14)$; 

\item [\rm (160)] $C_{32} (3,4,8,10,12,13,16)$, $C_{32} (4,5,8,10,11,12,16)$; 

\item [\rm (161)] $C_{32} (3,4,8,10,13,14,16)$, $C_{32} (4,5,8,10,11,14,16)$; 

\item [\rm (162)] $C_{32} (3,4,8,12,13,14,16)$, $C_{32} (4,5,8,11,12,14,16)$; 

\item [\rm (163)] $C_{32} (3,4,10,12,13,14,16)$, $C_{32} (4,5,10,11,12,14,16)$; 

\item [\rm (164)] $C_{32} (3,6,8,10,12,13,14)$, $C_{32} (5,6,8,10,11,12,14)$; 

\item [\rm (165)] $C_{32} (3,6,8,10,13,14,16)$, $C_{32} (5,6,8,10,11,14,16)$; 

\item [\rm (166)] $C_{32} (3,6,8,12,13,14,16)$, $C_{32} (5,6,8,11,12,14,16)$; 

\item [\rm (167)] $C_{32} (3,6,10,12,13,14,16)$, $C_{32} (5,6,10,11,12,14,16)$;

\item [\rm (168)] $C_{32} (3,8,10,12,13,14,16)$, $C_{32} (7,8,10,11,12,14,16)$; 
 
\item [\rm (169)] $C_{32} (2,3,4,6,8,10,12,13)$, $C_{32} (2,4,5,6,8,10,11,12)$; 

\item [\rm (170)] $C_{32} (2,3,4,6,8,10,13,16)$, $C_{32} (2,4,5,6,8,10,11,16)$;

\item [\rm (171)] $C_{32} (2,3,4,6,8,12,13,14)$, $C_{32} (2,4,5,6,8,11,12,14)$; 

\item [\rm (172)] $C_{32} (2,3,4,6,8,12,13,16)$, $C_{32} (2,4,5,6,8,11,12,16)$; 

\item [\rm (173)] $C_{32} (2,3,4,6,8,13,14,16)$, $C_{32} (2,4,5,6,8,11,14,16)$; 

 \item [\rm (174)] $C_{32} (2,3,4,6,10,12,13,16)$, $C_{32} (2,4,5,6,10,11,12,16)$; 

 \item [\rm (175)] $C_{32} (2,3,4,6,12,13,14,16)$, $C_{32} (2,4,5,6,11,12,14,16)$; 

 \item [\rm (176)] $C_{32} (2,3,4,8,10,12,13,14)$, $C_{32} (2,4,5,8,10,11,12,14)$; 

\item [\rm (177)] $C_{32} (2,3,4,8,10,12,13,16)$, $C_{32} (2,4,5,8,10,11,12,16)$; 

\item [\rm (178)] $C_{32} (2,3,4,8,10,13,14,16)$, $C_{32} (2,4,5,8,10,11,14,16)$; 

\item [\rm (179)] $C_{32} (2,3,4,10,12,13,14,16)$, $C_{32} (2,4,5,10,11,12,14,16)$; 

\item [\rm (180)] $C_{32} (2,3,6,8,10,12,13,16)$, $C_{32} (2,5,6,8,10,11,12,16)$; 

\item [\rm (181)] $C_{32} (2,3,6,8,12,13,14,16)$, $C_{32} (2,5,6,8,11,12,14,16)$; 

\item [\rm (182)] $C_{32} (2,3,8,10,12,13,14,16)$, $C_{32} (2,5,8,10,11,12,14,16)$; 

\item [\rm (183)] $C_{32} (3,4,6,8,10,12,13,14)$, $C_{32} (4,5,6,8,10,11,12,14)$; 

\item [\rm (184)] $C_{32} (3,4,6,8,10,13,14,16)$, $C_{32} (4,5,6,8,10,11,14,16)$; 

\item [\rm (185)] $C_{32} (3,4,6,8,12,13,14,16)$, $C_{32} (4,5,6,8,11,12,14,16)$; 

\item [\rm (186)] $C_{32} (3,4,6,10,12,13,14,16)$, $C_{32} (4,5,6,10,11,12,14,16)$; 

\item [\rm (187)] $C_{32} (3,4,8,10,12,13,14,16)$, $C_{32} (4,5,8,10,11,12,14,16)$; 

\item [\rm (188)] $C_{32} (3,6,8,10,12,13,14,16)$, $C_{32} (5,6,8,10,11,12,14,16)$; 

\item [\rm (189)] $C_{32} (2,3,4,6,8,10,12,13,16)$, $C_{32} (2,4,5,6,8,10,11,12,16)$; 

\item [\rm (190)] $C_{32} (2,3,4,6,8,12,13,14,16)$, $C_{32} (2,4,5,6,8,11,12,14,16)$; 

\item [\rm (191)] $C_{32} (2,3,4,8,10,12,13,14,16)$, $C_{32} (2,4,5,8,10,11,12,14,16)$; 

\item [\rm (192)] $C_{32} (3,4,6,8,10,12,13,14,16)$, $C_{32} (4,5,6,8,10,11,12,14,16)$; 

\end{enumerate}

\section{Conclusion}

 The author feels that this paper provides a lot of scope for further research to obtain more families of Type-2 isomorphic circulant graphs.
 
\vspace{.1cm}
\noindent
\textbf{Declaration of competing interest}\quad 
The author declares that he has no conflict of interest.


\end{document}